\documentclass[letterpaper,11pt,oneside,reqno]{amsart}
\usepackage{bm}
\usepackage{mathrsfs}
\usepackage{amsfonts,amsmath, amssymb,amsthm,amscd,stmaryrd}
\usepackage[height=9.6in,width=5.95in]{geometry}
\usepackage[mpexclude,DIV13]{typearea}
\usepackage{verbatim}
\usepackage{hyperref}
\usepackage{graphicx}
\usepackage[latin1]{inputenc}
\usepackage{latexsym}
\usepackage{lscape}
\usepackage{epsfig}
\include{bibtex}
\usepackage[usenames,dvipsnames]{color}

\usepackage{setspace}

\usepackage{parskip}
\newtheorem{theorem}{Theorem}[section]
\newtheorem{lemma}[theorem]{Lemma}

\newtheorem{corollary}[theorem]{Corollary}
\newtheorem{proposition}[theorem]{Proposition}

\newtheorem{definition}[theorem]{Definition}

\def\N{\mathbb{N}}

\def\P{\mathbb{P}}
\def\Z{\mathbb{Z}}
\def\R{\mathbb{R}}

\def\E{\mathbb{E}}

\newcommand{\PP}{\ensuremath{\mathbb{P}}}

\newcommand{\wgt}[3]{W_{#1}(#2,#3)}
\newcommand{\polymer}[3]{\rho_{#1}(#2,#3)}
\newcommand{\Cpop}{G}
\newcommand{\maxpoly}{\mathrm{MaxDisjtPoly}}
\newcommand{\maxpolymac}[3]{{\mathrm{MaxDisjtPoly}}_{#1}(#2,#3)}

\newcommand{\mc}{\mathcal}

\newcommand{\e}{\varepsilon}

\begin{document}
\title[Fractal geometry of the Airy sheet]{Fractal geometry of Airy$_2$ processes coupled~via~the~Airy~sheet}

\author[R. Basu]{Riddhipratim Basu}
\address{R. Basu\\
  International Centre for Theoretical Sciences\\
  Tata Institute of Fundamental Research\\
  Bangalore, India
  }
  \email{rbasu@icts.res.in}

\author[S. Ganguly]{Shirshendu Ganguly}
\address{S. Ganguly\\
  Department of Statistics\\
 U.C. Berkeley \\
  Evans Hall \\
  Berkeley, CA, 94720-3840 \\
  U.S.A.}
  \email{sganguly@berkeley.edu}

\author[A. Hammond]{Alan Hammond}
\address{A. Hammond\\
  Department of Mathematics and Statistics\\
 U.C. Berkeley \\
  Evans Hall \\
  Berkeley, CA, 94720-3840 \\
  U.S.A.}
  \email{alanmh@berkeley.edu}

\begin{abstract}
In last passage percolation models lying in the Kardar-Parisi-Zhang universality class, maximizing paths that travel over distances of order $n$ accrue energy that fluctuates on scale $n^{1/3}$; and these paths deviate from the linear interpolation of their endpoints on scale $n^{2/3}$. 
These maximizing paths and their energies may be viewed via a coordinate system that respects these scalings. What emerges by doing so is a system  indexed by $x,y \in \R$ and $s,t \in \R$ with $s < t$  of unit order quantities $W_n\big( x,s ; y,t \big)$ specifying the scaled energy of the maximizing path that moves in scaled coordinates between $(x,s)$ and $(y,t)$. 
The space-time Airy sheet   is, after a parabolic adjustment,  the putative distributional limit $W_\infty$ of this system as $n \to \infty$. The Airy sheet has recently been constructed in~\cite{DOV18}  as such a limit of Brownian last passage percolation. In this article, we initiate the study of fractal geometry in the Airy sheet. We prove that the scaled energy difference profile given by  $\R \to \R: z \to W_\infty \big( 1,0 ; z,1 \big) - W_\infty \big( -1,0 ; z,1 \big)$ is a non-decreasing process that is constant in a random neighbourhood of almost every $z \in \R$; and that the exceptional set of $z \in \R$ that violate this condition almost surely has Hausdorff dimension one-half. Points of violation correspond to special behaviour for scaled maximizing paths, and we prove the result by investigating this behaviour, making use of two inputs from recent studies of scaled Brownian LPP; namely,  Brownian regularity of profiles, and estimates on the rarity of pairs of disjoint scaled maximizing paths that begin and end close to each other.  
\end{abstract}

\maketitle
\tableofcontents

\section{Introduction}
\subsection{Kardar-Parisi-Zhang universality and last passage percolation}

The Kardar-Parisi-Zhang [KPZ] equation is a stochastic PDE putatively modelling 
 a wide array of models of one-dimensional local random growth subject to restraining forces such as surface tension. The theory of KPZ universality 
 predicts that these models share
a triple $(1,1/3,2/3)$ of exponents: in time of scale $t^1$, a growing interface  above a given point in its domain differs from its mean value by a height that is a random quantity of order $t^{1/3}$; and it is by varying this point on a spatial scale of $t^{2/3}$ that non-trivial correlation between the associated random heights is achieved. When the random height over a given point is scaled by dividing by $t^{1/3}$, a scaled quantity is obtained whose limiting law in high $t$  
is governed by the extreme statistics of certain ensembles of large random matrices. The theory has been the object of several intense waves of mathematical attention in recent years. One-parameter models whose parameter corresponds to time in the KPZ equation have been rigorously analysed~\cite{ACQ11,OC12,COSZ14,BC14} by integrable techniques so that the equation is seen to describe the evolution of certain weakly asymmetric growth models. Hairer's theory of regularity structures~\cite{Hairer14} has provided a robust concept of solution to the equation, raising (and realizing:~\cite{HQ15}) the prospect of deriving invariance principles.

In 
last passage percolation [LPP] models, a random environment which is independent in disjoint regions is used to assign random values called energies to paths that run through it. A path with given endpoints of maximal energy is called a geodesic. LPP is concerned with the behaviour of energy and geometry of geodesics that run between distant endpoints. The large-scale behaviour of many LPP models is expected to be governed by the KPZ exponent triple -- pioneering rigorous works concerning Poissonian LPP are~\cite{BDJ99} and \cite{J00} -- and it is natural to view these models through the lens of scaled coordinates 
whose choice is dictated by this triple. Since LPP models thus scaled are expected to be described by a scaled form of the KPZ equation in the limit of late time, they offer a suitable framework for the study of the KPZ fixed point, namely of those random objects which are shared between models in the KPZ universality class by offering an accurate scaled description of large-scale and late-time behaviour of such models. 

The study of last passage percolation in scaled coordinates depends critically on inputs of integrable origin, but it has been recently proved profitable to advance it through several  probabilistic perspectives on KPZ universality. It will become easier to offer signposts to pertinent articles after we have specified the LPP model that we will study.
This paper is devoted to giving rigorous expression to a novel aspect of the KPZ fixed point: to the fractal geometry of the stochastic process given by the difference in scaled energy of a pair of geodesics rooted at given fixed distant horizontally displaced lower endpoints as the higher endpoint, which is shared between the two geodesics, is varied horizontally. The concerned result is proved by exploiting and developing recent advances in the rigorous theory of KPZ in which probabilistic tools are harnessed in unison with limited but essential integrable inputs.  

We next present the Brownian last passage percolation model that will be our object of study; explain how it may be represented in scaled coordinates; briefly discuss recent probabilistic tools in KPZ; and state our main theorem.

\subsection{Brownian last passage percolation: geodesics and their energy}\label{s.brlpp}
In this LPP model, a field of local randomness is specified by 
 an ensemble $B:\Z \times \R \to \R$ of independent  two-sided standard Brownian motions $B(k,\cdot):\R\to \R$, $k \in \Z$, defined 
on a probability space carrying a law that we will label~$\PP$. 

Any non-decreasing path $\phi$ mapping a compact real interval to $\Z$ is ascribed an energy~$E(\phi)$ by summing the Brownian increments associated to~$\phi$'s level sets.
To wit,
let $i,j \in \Z$ with $i \leq j$.
We denote the integer interval $\{i,\cdots,j\}$ by $\llbracket i,j \rrbracket$.
Further let $x,y \in \R$ with $x \leq y$.
Each non-decreasing function $\phi: [x,y] \to \llbracket i,j \rrbracket$
with $\phi(x) = i$ and $\phi(y) = j$
corresponds to a non-decreasing list 
 $\big\{ z_k: k \in \llbracket i+1,j \rrbracket \big\}$ of values $z_k \in [x,y]$ if we select $z_k = \sup \big\{ z \in [a,b]: \phi(z) \leq k-1 \big\}$. 
With the convention that $z_i = x$ and $z_{j+1} = y$,
the path energy $E(\phi)$ is set equal to $\sum_{k=i}^j \big( B ( k,z_{k+1} ) - B( k,z_k ) \big)$.
We then define  the maximum energy
$$
M_{(x,i) \to (y,j)} \, = \, \sup \, \bigg\{ \, \sum_{k=i}^j \Big( B ( k,z_{k+1} ) - B( k,z_k ) \Big) \, \bigg\} \, , 
$$
where this supremum of energies $E(\phi)$ is taken over all such paths~$\phi$. The random process $M_{(0,1) \to (\cdot,n)}: [0,\infty) \to \R$ was introduced by~\cite{GlynnWhitt} and further studied in~~\cite{baryshnikov2001gues} and~\cite{OCY}.

It is perhaps useful to visualise a non-decreasing path such as $\phi$ above by viewing it as the associated {\em staircase}. The staircase associated to $\phi$ is a subset of the planar rectangle $[x,y] \times [i,j]$ given by the range of a continuous path between $(x,i)$ and $(y,j)$ that alternately moves horizontally and vertically. The staircase is a union of horizontal and vertical planar line segments. The horizontal segments are $[z_k,z_{k+1}] \times \{ k \}$ for $k \in \llbracket i , j\rrbracket$; while the vertical segments interpolate the right and left endpoints of consecutively indexed horizontal segments. 


\subsection{Scaled coordinates for Brownian LPP: polymers and their weights}

The one-third and two-thirds KPZ scaling considerations are manifest in Brownian LPP. When the ending height~$j$ exceeds the starting height $i$ by a large quantity $n \in \N$, and the location $y$ exceeds $x$ also by $n$, then the maximum energy grows linearly, at rate $2n$,
and has a fluctuation about this mean of order~$n^{1/3}$. 
Indeed, the maximum energy of any path of journey $(0,0) \to (n,n)$ 
verifies 
\begin{equation}\label{e.energyscaled}
M_{(0,0) \to (n,n)}  = 2n +  2^{1/2} n^{1/3} W_n \, ,
\end{equation}
where $W_n$ is a scaled expression for energy; since $M_{(0,0) \to (n,n)}$  has the law of
the uppermost particle at time $n$ in a Dyson Brownian motion with $n+1$ particles by~\cite[Theorem~$7$]{OCY}, and the latter law has the distribution of the uppermost eigenvalue of an $(n+1) \times (n+1)$ matrix randomly drawn from the Gaussian unitary ensemble with entry variance $n$ by~\cite[Theorem~$3$]{Grabiner}, the quantity $W_n$ converges in distribution as $n \to \infty$
to the Tracy-Widom GUE distribution. 
Any non-decreasing path $\phi: [0,n] \to \llbracket 0 , n \rrbracket$ that attains this maximal energy will be called a geodesic, and denoted for use in a moment by~$P_n$; the term geodesic is further applied to any non-decreasing path that realizes the maximal energy assumed by such paths that share its initial and final values. 

Moreover, when the horizontal coordinate of the ending point of the journey $(0,0) \to (n+y,n)$ is permitted to vary away from $y=0$, then it is changes of $n^{2/3}$ in the value of~$y$ that result in a non-trivial correlation of the maximum energy with its original value.

Universal large-scale properties of LPP may be studied by using scaled coordinates to depict geodesics and their energy; a geodesic thus scaled will be called a {\em polymer} and its scaled energy will be called its {\em weight}.

\begin{figure}
\includegraphics[width=0.7\textwidth]{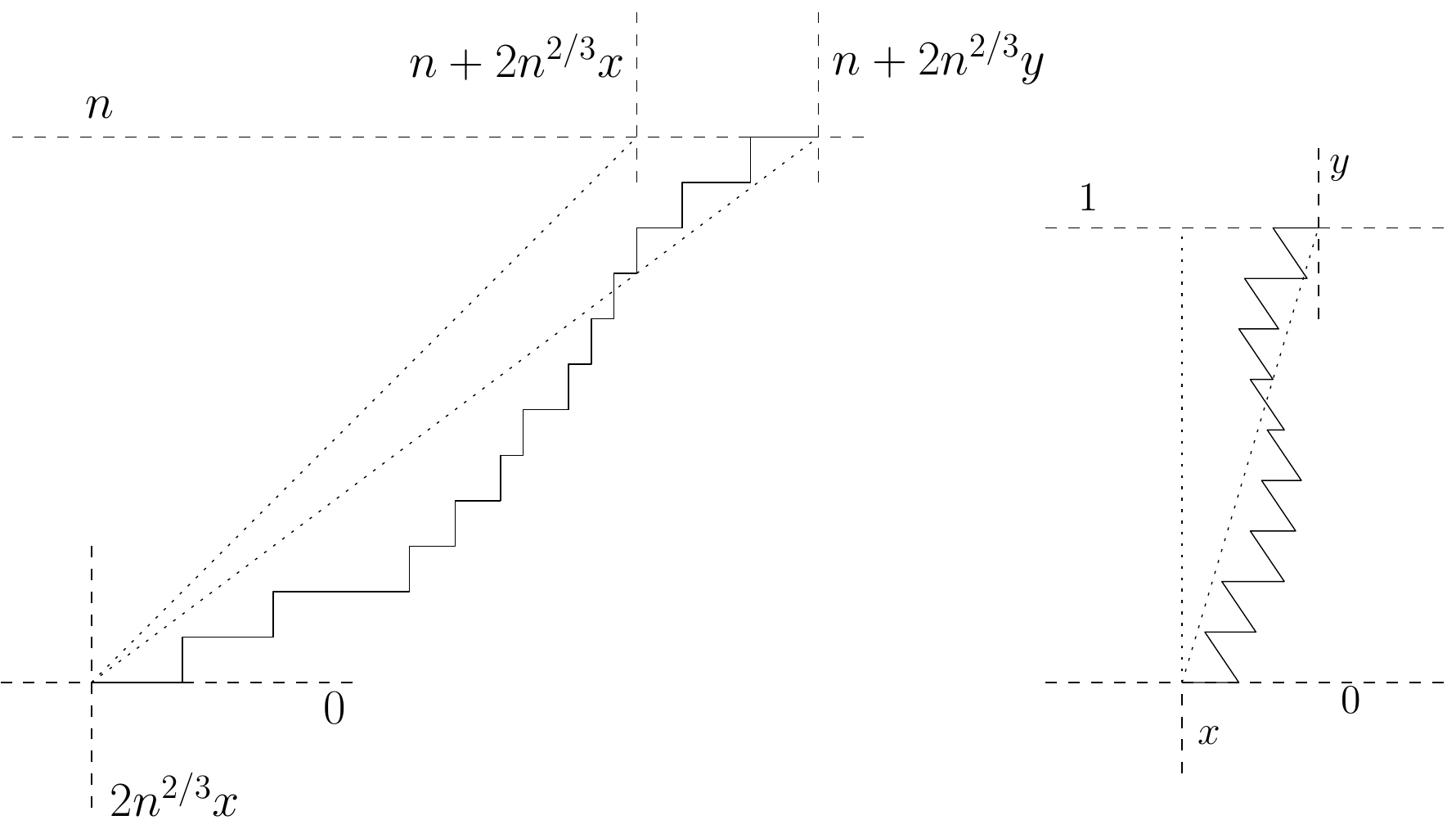}
\caption{The staircase in the left sketch is associated to a geodesic. When the  scaling map $R_n$ is applied to it, the outcome is the polymer $\rho_n(x,y)$ in the right sketch.}\label{f.scaling}
\end{figure}

Let $R_n:\R^2 \to \R^2$ be the {\em scaling map}, namely the linear map sending $(n,n)$ to $(0,1)$ and $(2n^{2/3},0)$ to $(1,0)$.  The image of any staircase under the scaling map will be called an $n$-zigzag. An $n$-zigzag is comprised of planar line segments that are consecutively  horizontal and downward sloping but near horizontal, the latter type each having gradient $-2n^{-1/3}$.
For $x,y \in \R$ and $n \in \N$, let $\rho_n(x,y)$, a subset of $\R \times [0,1]$,
denote the image under $R_n$
of the staircase associated to the LPP geodesic whose endpoints are 
$R_n^{-1} (x,0)$ and $R_n^{-1} (y,1)$. For example, $\rho_n(0,0)$ is the image under the scaling map of the staircase attached to the geodesic $P_n$; for  $x,y \in \R$, $\rho_n(x,y)$ is the $n$-polymer (or scaled geodesic) which crosses the unit-strip $\R \times [0,1]$ between $(x,0)$ and $(y,1)$. For any given pair $(x,y)$ that we consider, $\rho_n(x,y)$ is well defined, because there almost surely exists a unique $n$-polymer from $(x,0)$ to $(y,1)$ by \cite[Lemma~$4.6(1)$]{hammond2017patchwork}. 
This $n$-polymer is depicted in Figure~\ref{f.scaling}. The label $n$ is used consistently when $n$-zigzags and $n$-polymers are considered, and we refer to them simply as zigzags and polymers.

Scaled geodesics have scaled energy: in the case of $\rho_n(0,0)$, its scaled energy is the quantity $W_n$ appearing in~(\ref{e.energyscaled}). We now set $W_n(0,0) := W_n$ with a view to generalizing. Indeed, for $x,y \in \R$ satisfying $y - x \geq -2^{-1}n^{1/3}$, the {\em unit-order} scaled energy or {\em weight} $W_n(x,y)$ of $\rho_n(x,y)$ is given by 
\begin{equation}\label{e.weight}
  W_n(x,y) = 2^{-1/2}n^{-1/3} \Big( M_{(2n^{2/3}x,0) \to (n + 2n^{2/3}y,n)} - 2n - 2n^{2/3}(y-x) \Big) \, , 
\end{equation}
where note that $M_{(2n^{2/3}x,0) \to (n + 2n^{2/3}y,n)}$ is equal to the energy of the LPP geodesic whose staircase maps to $\rho_n(x,y)$ under~$T_n$.

The random weight profile of scaled geodesics emerging from $(0,0)$ to reach the horizontal line at height one, namely
$\R \to \R: z \to W_n (0,z)$, converges -- see \cite[Theorem~$1.2$]{J03} for the case of geometric LPP 
-- in the high $n$ limit to a canonical object in the theory of KPZ universality. This object, which is the Airy$_2$ process after the subtraction of a parabola $x^2$, has finite-dimensional distributions specified by Fredholm determinants. (It is in fact incorrect to view the domain of such profiles as $W_n(0,\cdot)$ as the whole of $\R$, but we tolerate this abuse until correcting it shortly.)

\subsection{Probabilistic and geometric approaches to last passage percolation}\label{s.probgeom}

This determinantal information about profiles such as 
$\R \to \R: z \to W_n (0,z)$ offers a rich store of exact formulas which nonetheless has not {\em per se} led to derivations of certain putative properties of this profile, such as Johansson's conjecture that its high~$n$ weak limit, the parabolically adjusted Airy$_2$ process, has an almost surely unique maximizer; or the absolute continuity, uniformly in high~$n$, of the profile on a compact interval with respect to a suitable vertical shift of Brownian motion. Probabilistic and geometric perspectives on LPP, allied with integrable inputs, have led to several recent advances, including the solution of these problems. The above profile may be embedded~\cite{baryshnikov2001gues,OCY} as the uppermost curve in an $n$-curve ordered ensemble of curves whose law is that of a system of Brownian motions conditioned on mutual avoidance subject to a suitable boundary condition. As such, this uppermost curve enjoys a simple {\em Brownian Gibbs} resampling property when it is resampled on a given compact interval in the presence of data from the remainder of the ensemble. The Brownian Gibbs property has been exploited in~\cite{CH14} to prove Brownian absolute continuity of the Airy$_2$ process as well as Johansson's conjecture -- and this conjecture has been obtained by   Moreno  Flores,  Quastel  and  Remenik~\cite{MFQR13}  via  an  explicit  formula  for
the maximizer, and an argument of Pimentel~\cite{Pimentel14} showing that any stationary process minus a
parabola has a unique maximizer.  A  positive temperature analogue of the Brownian Gibbs property has treated~\cite{CH14} questions of Brownian similarity for the scaled narrow wedge solution of the KPZ equation; and a more refined understanding~\cite{H16} of Brownian regularity of the profile  $\R \to \R: z \to W_n (0,z)$  has been obtained by further Brownian Gibbs analysis.

Robust probabilistic tools harnessing merely integrable one-point tail bounds have been used to study non-integrable perturbations of LPP problems, such as in the solution~\cite{BSS14} of the slow bond problem; bounds on coalescence times for LPP geodesics~\cite{BSS17++}; and to identify~\cite{FerrariSpohn,FO18,BG18} a Holder exponent of $1/3-$ for the weight profile when the latter endpoint is varied in the vertical, or temporal, direction. 
As~\cite{Sep18} surveys,  geometric properties such as fluctuation and coalescence of geodesics have been studied~\cite{BCS06, Pim16} in stationary versions of LPP by using queueing theory and the Burke property.

The distributional convergence in high $n$ of the profile $\R \to \R: z \to W_n (0,z)$ --  and counterpart convergences for certain other integrable LPP models -- to a limiting stochastic process is by now a classical part of the rigorous theory of LPP. It has  expected since at least~\cite{CQR15} that a richer universality object, the {\em space-time Airy sheet}, specifying the limiting weight of polymers between pairs of  planar points $(x,u)$ and $(y,v)$ that are arbitrary except for the condition that $u \not= v$, should exist uniquely. Two significant recent advances address this and related universal objects. The polymer weight profile $\R \to \R: z \to W_\infty (0,z)$ may be viewed as the limiting time-one snapshot of an evolution in positive time begun from the special initial condition consisting of a Dirac delta mass at the origin. In the first advance~\cite{MQR17}, this evolution is constructed for all positive time  from an almost arbitrary general initial condition (in fact, the totally asymmetric exclusion process is used as the prelimiting model, in place of Brownian LPP); explicit Fredholm determinant formulas for the evolution are provided. (The Brownian regularity of the time-one snapshot of this evolution from general initial data is studied in~\cite{hammond2017patchwork} for the Brownian LPP prelimit.)
In the second recent contribution~\cite{DOV18}, the entire space-time Airy sheet is constructed, by use of an extension of the Robinson-Schensted-Knuth correspondence which permits the sheet's construction in terms of a last passage percolation problem whose underlying environment is itself a copy of the high~$n$ distributional limit of the narrow wedge profile $\R \to \R: z \to W_n (0,z)$. The analysis of~\cite{DOV18} is assisted by~\cite{DV18}, an article making a Brownian Gibbs analysis of scaled Brownian  LPP in order to provide valuable estimates for the study of the very novel LPP problem introduced in~\cite{DOV18}.   

\subsection{The main result, concerning fractal random geometry in scaled Brownian LPP}


All of the above is to say that robust probabilistic tools have furnished a very fruitful arena in the recent study of scaled LPP problems. In the present article, we isolate an aspect of the newly constructed Airy sheet in order to shed light on the fractal geometry of this rich universal object. We will use the lens of the prelimiting scaled Brownian LPP model to express our principal result, and then record a corollary that asserts the corresponding statement about fractal geometry in the Airy sheet.

\begin{figure}
\includegraphics[width=0.7\textwidth]{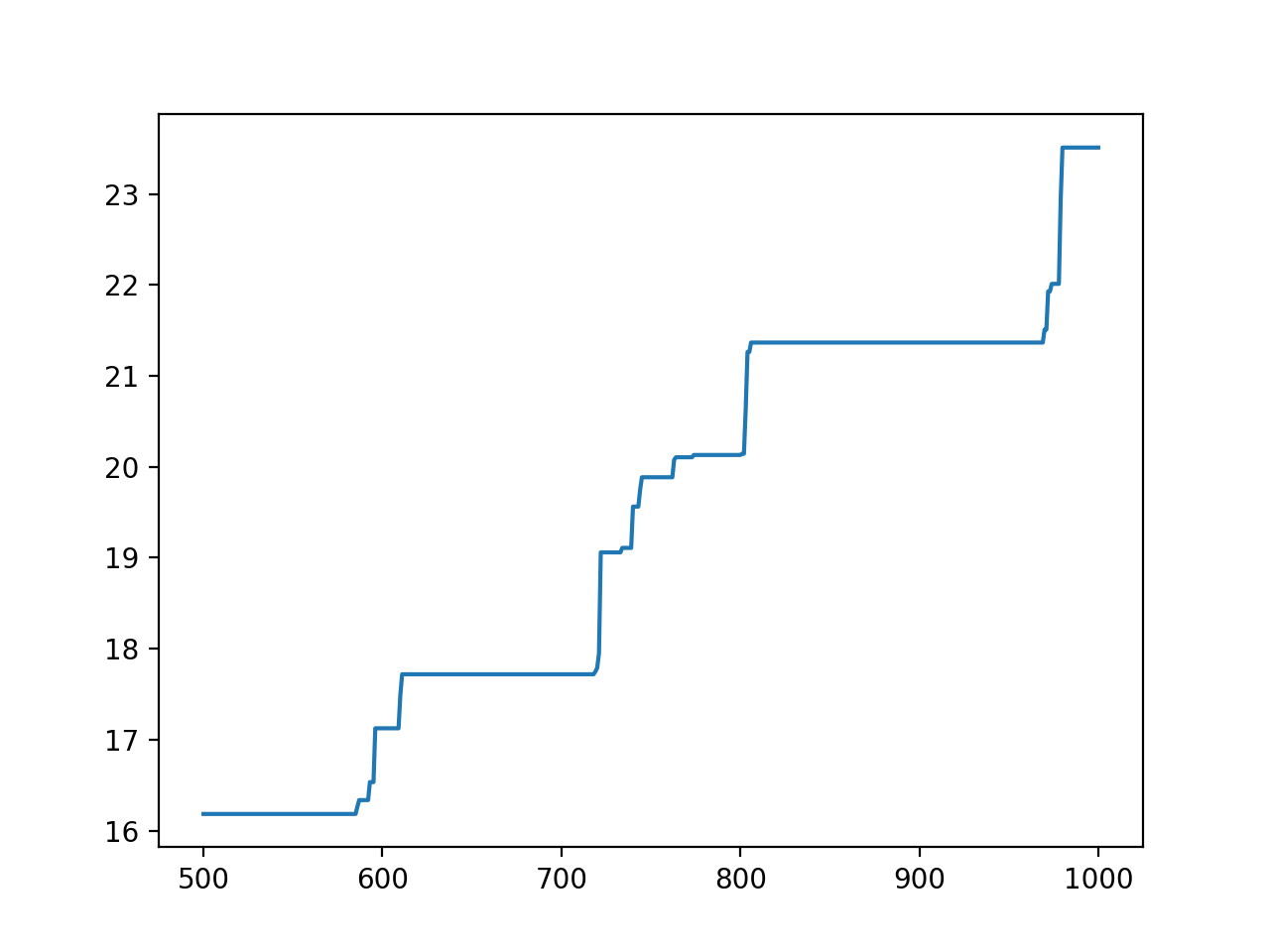}
\caption{
A simulation of an LPP model by 
Junou Cui,  Zoe Edelson and Bijan Fard. With $n = 500$, the difference in energy of geodesics making the unscaled journeys $(-n^{2/3},0) \to (z,n)$ and $(n^{2/3},0) \to (z,n)$ is plotted against $x$-coordinate~$z$.}\label{f.simulation} 
\end{figure}

The novel process that is our object of study is the random {\em difference weight profile} given by considering the relative weight of unit-height polymers in Brownian LPP emerging from the points $(-1,0)$ and $(1,0)$; namely, $z  \to W_n(1,z) - W_n(-1,z)$. This real-valued stochastic process is defined under the condition $z \geq -2^{-1}n^{1/3} + 1$ that ensures that the constituent weights are well specified by~(\ref{e.weight}); but we may extend the process' domain of definition to the whole of the real line by setting it equal to its value at $-2^{-1}n^{1/3} + 1$ for smaller $z$-values. Since the random functions $z \to W_n(x,z)$ for each $x \in \R$ are almost surely continuous by~\cite[Lemma 2.2(1)]{hammond2017modulus}, we see that $\R \to \R: z \to W_n(1,z) - W_n(-1,z)$ is almost surely an element of the space $\mc{C}$ of real-valued continuous functions on $\R$. Equipping $\mc{C}$ with the topology of locally uniform convergence, we may consider weak limit points in the limit of high $n$ of this difference weight profile.  
Our principal result asserts that such limit points are the distribution functions of random Cantor sets: see Figure~\ref{f.simulation} for a simulation.
\begin{theorem}\label{t.z}
Any weak limit as $n \to \infty$ of the sequence of random processes $\R \to \R: z \to W_n(1,z) - W_n(-1,z)$ is a random function $Z:\R \to \R$
such that
\begin{enumerate} 
\item $Z$ is almost surely continuous and non-decreasing;
\item $Z$ is constant in a random neighbourhood of almost every 
$z \in \R$; 
\item the set $E$ of points $z \in \R$ 
that violate the preceding condition -- those $z$ about which $Z$ is not locally constant --
is thus a Lebesgue null set a.s.; this set almost surely has Hausdorff dimension one-half. 
\end{enumerate}
\end{theorem}

We have expressed our principal result in the language of weak limit points in order to permit its adaptation to other LPP models; and because it is natural to prove the result by deriving counterpart assertions (which may be useful elsewhere) for the Brownian LPP prelimit. 
Significantly, however, the Airy sheet has been constructed; moreover, it is Brownian LPP which is the prelimiting model in that construction. We may thus express a corollary in terms of the Airy sheet. 

That is, 
let $W = W_\infty: \R^2 \to \R$ denote the parabolically shifted Airy sheet constructed in~\cite{DOV18}. Namely, endowing the space of continuous real-valued functions on $\R^2$ with the topology of locally uniform convergence, the random function $(x,y) \to W_n(x,y)$ converges weakly to $W:\R^2 \to \R$ as $n \to \infty$ 
by \cite[Theorem~$1.3$]{DOV18}. 
\begin{corollary} 
Theorem~\ref{t.z} remains valid when the random function $Z: \R \to \R$
is set equal to $Z(z) = W(1,z) - W(-1,z)$. 
\end{corollary}
{\bf Proof.} Since the weak limit point $Z$ in Theorem~\ref{t.z} exists, is unique, and is given in law  by $\R \to \R: z \to W(1,z) - W(-1,z)$, the result follows from the theorem. \qed

Theorem~\ref{t.z} concerns the fractal geometry of random Cantor sets that are  embedded in a canonical universal object that arises 
as a scaling limit of statistical mechanical models. It shares these features with the distribution function of the local time at zero of one-dimensional Brownian motion -- in fact, as~\cite[Theorem~$4.24$]{PeresMorters} shows, even the one-half Hausdorff dimension of the random Cantor set is shared with this simple example. The qualitative features are also shared with the distribution function associated to a natural local time constructed~\cite{HPS} on the set of exceptional times of dynamical critical percolation on the hexagonal lattice, in which case, the Hausdorff dimension of the exceptional set is known by~\cite{GPS} to equal $31/36$.  

Although we have presented Theorem~\ref{t.z} as our main result, its proof will yield an interesting consequence, establishing the sharpness of the exponent in a recent upper bound on the probability of the presence of a pair of disjoint polymers that begin and end at nearby locations. Since we will anyway review the concerned upper bound in the next section, we defer the statement of this second theorem to Section~\ref{s.two}.

\subsection{Acknowledgments}
The authors thank Duncan Dauvergne, Milind Hegde and B{\'a}lint Vir{\'a}g for helpful comments concerning a draft version of this paper.
R.B. is partially supported by a Ramanujan Fellowship from the Government of India, and an ICTS-Simons Junior Faculty Fellowship. A.H. is supported by National Science Foundation grant DMS-$1512908$ and by a Miller Professorship. This work was conducted in part during a visit of R.B. to the Statistics Department at U.C. Berkeley, whose hospitality he gratefully acknowledges.

\section{Brownian weight profiles and disjoint polymer rarity; \\ and a~conceptual~overview~of~the~main~argument}\label{s.two}

In the first two subsections, we provide the two principal inputs for our main result; in a third, we explain in outline how to use them to prove it; in a fourth, we state our second principal result, Theorem~\ref{t.disjtpoly.lb}; and, in the fifth, we record some basic facts about polymers.

\subsection{Polymer weight change under horizontal perturbation of endpoints}\label{s.pwc}
Set $Q:\R \to \R$ equal to the parabola $Q(z) = 2^{-1/2} z^2$.
For any given $x \in \R$, the polymer weight profile   $y \to \wgt{n}{x}{y}$  has in the large scale a curved shape that in an average sense peaks at $x$, the profile hewing to the curve $-Q(y-x)$. When this parabolic term is added to the polymer weight, the result is a random process in $(x,y)$
which typically suffers changes of order $\e^{1/2}$
when $x$ or $y$ are varied on a small scale $\e > 0$.
Our first main input gives rigorous expression to this statement, uniformly in $(n,x,y) \in \N \times \R \times \R$ for which the difference $\vert y - x \vert$ is permitted to inhabit an expanding region about the origin, of scale~$n^{1/18}$.

\begin{theorem}\cite[Theorem~$1.1$]{hammond2017modulus}\label{t.differenceweight}
Let $\e \in (0,2^{-4}]$. 
Let $n \in \N$ satisfy
$n \geq 10^{32} c^{-18}$, and let $x,y \in \R$ verify  $\big\vert x - y  \big\vert \leq 2^{-2} 3^{-1} c  n^{1/18}$.
Let 
 $R \in \big[10^4 \, , \,   10^3 n^{1/18} \big]$.
Then
$$
\PP \left( \sup_{\begin{subarray}{c} u \in [x,x+\e] \\
    v \in [y,y+\e]  \end{subarray}} \Big\vert \wgt{n}{u}{v} + Q(v - u) - \wgt{n}{x}{y} - Q(y - x) \Big\vert  \, \geq \, \e^{1/2}
  R  \right)
$$
  is at most  $10032 \, C  \exp \big\{ - c_1 2^{-21}   R^{3/2}   \big\}$, where  $c_1 = \min \big\{ 2^{-5/2} c , 1/8 \big\}$.
\end{theorem}
The bound in Theorem~\ref{t.differenceweight}, and several later results, have been expressed explicitly up to two positive constants $c$ and $C$. We reserve these two symbols for this usage throughout. 
The concerned pair of constants enter via bounds that we will later quote in Theorem~\ref{t.onepoint}  on the upper and lower tail of the uppermost eigenvalue of an $n \times n$ matrix randomly selected from the Gaussian unitary ensemble.

The imposition in Theorem~\ref{t.differenceweight} that 
 $R \leq 10^3 n^{1/18}$ is rather weak in the sense much of the result's interest lies in high choices of $n$. Indeed, we now provide a formulation in which this condition is absent. 

\begin{corollary}\label{c.aest} 
There exist positive constants $C_1$ and $c_2$ such that, for $\alpha \in (0,1/2)$ and  $\e \in (0,2^{-4}]$,
\begin{equation}
 \limsup_{n \in \N} \, \P \left( \sup_{\begin{subarray}{c}
    v_1,v_2 \in [y,y+\e]  \end{subarray}} \Big\vert \wgt{n}{0}{v_2}+Q(v_2)-  \wgt{n}{0}{v_1}-Q(v_1) \Big\vert  \geq \e^\alpha \right) \leq C_1 \exp  \big\{ - c_2 \e^{3(2\alpha-1)/4} \big\} \, .
\end{equation}
\end{corollary}
{\bf Proof.}
Set $c_2$ equal to the quantity $2^{-21} c_1$ from Theorem~\ref{t.differenceweight}. Then apply this result with 
 $R$ set equal to $\e^{(2\alpha-1)/2}$, choosing $C_1$ high enough that the hypothesis $R \geq 10^4$ may be supposed due to the desired result being vacuously satisfied in the opposing case. \qed

\subsection{The rarity of pairs of polymers with close endpoints}

Let $n \in \N$ and let $I,J \subset \R$ be intervals. Set $\maxpolymac{n}{I}{J}$ equal to  the maximum cardinality of a {\em pairwise disjoint} set of   polymers each of whose starting and ending points have the respective forms $(x,0)$ and $(y,1)$ 
where $x$ is some element of $I$ and $y$ is some element of $J$.

The second principal input gauges the rarity of the event that this maximum cardinality exceeds any given $k \in \N$ when $I$ and $J$ have a given short length $\e$. We will apply the input with $k=2$, since it is the rarity of pairs of polymers with nearby starting and ending points that will concern~us.

\begin{theorem}\cite[Theorem~$1.1$]{brownianLPPtransversal}\label{t.disjtpoly.pop}
There exists a positive  constant $\Cpop$ such that the following holds.
Let  $n \in \N$; and 
let  $k \in \N$, $\e > 0$ and $x,y \in \R$ satisfy the conditions that $k \geq 2$,
$$
 \e \leq    \Cpop^{-4k^2}   \, , \, n  \geq   \Cpop^{k^2}  \big( 1 +   \vert x - y \vert^{36}    \big) \e^{-\Cpop} 
$$
and  
$\vert x - y \vert  \leq \e^{-1/2} \big( \log \e^{-1} \big)^{-2/3} G^{-k}$.

Setting $I = [x- \e,x+ \e]$ and $J = [y- \e,y+ \e]$, we have that  
$$
\PP \Big( \maxpolymac{n}{I}{J}  \geq k \Big) 
  \leq  
 \e^{(k^2 - 1)/2} \cdot R \, ,
$$
where $R$ is a positive correction term that is bounded above by $\Cpop^{k^3}  \exp \big\{ \Cpop^k \big( \log \e^{-1} \big)^{5/6} \big\}$.
\end{theorem}
An alternative regime, where $\e$ is of unit order and $k \in \N$ is large, is addressed by~\cite[Theorem~$2$]{BHS18}: the counterpart of $\PP \big( \maxpolymac{n}{[-1,1]}{[-1,1]} \geq k \big)$ is bounded above by $\exp \big\{ - d k^{1/4} \big\}$ for some positive constant~$d$. 

Both inputs Theorem~\ref{t.differenceweight} and Theorem~\ref{t.disjtpoly.pop} have derivations
depending on the Brownian Gibbs resampling technique that we mentioned in Subsection~\ref{s.probgeom}. This use is perhaps more fundamental in the case of Theorem~\ref{t.disjtpoly.pop}, whose proof operates by showing that the presence of $k$ polymers with $\e$-close endpoint pairs typically entails a near touch of closeness of order~$\e^{1/2}$ at a given point on the part of the uppermost $k$ curves in the ordered ensembles of curves to which we alluded in Subsection~\ref{s.probgeom}; 
Brownian Gibbs arguments provide an upper bound on the latter event's probability.

\subsection{A conceptual outline of the main proof}

Theorem~\ref{t.z} is proved by invoking the two results just cited. Here we explain roughly how, thus explicating how our result on fractal geometry in scaled LPP is part of an ongoing probabilistic examination of universal KPZ objects. 

The theorem has three parts, and our heuristic discussion of the result's proof treats each of these in turn.

\subsubsection{Heuristics $1$: continuity of the weight difference profile.} The limiting profile $Z$ is a difference of parabolically shifted Airy$_2$ processes (which are coupled together in a non-trivial way). Since the Airy$_2$ process is almost surely continuous, so is $Z$. That $Z$ is non-decreasing is a consequence of a short planarity argument of which we do not attempt an overview, but  which has appeared in the proof of~\cite[Proposition~$3.8$]{DOV18}; and a variant of which originally addressed problems in first passage percolation~\cite{AlmWierman}.

\subsubsection{Heuristics $2$: local constancy of $Q$ about almost every point.} Logically, Theorem~\ref{t.z}(2) is merely a consequence of Theorem~\ref{t.z}(3), but it may be helpful to offer a guide to a proof in any case.
Implicated in the assertion is the geometric behaviour of the associated pair of random fields of polymers
$\big\{ \rho_n(\pm 1,z): z \in \R \big\}$.
For given $z \in \R$, $\rho_n(-1,z)$ and $\rho_n(1,z)$
respectively leave $(-1,0)$ and $(1,0)$. They arrive together at $(z,1)$
having merged at some random intermediate height $h_n \in (0,1)$. A key {\em coalescence} observation -- 
which we will not verify directly in the actual proof but which is a close cousin of Theorem~\ref{t.disjtpoly.pop}
-- is that this merging occurs in a uniform sense in $n$ away from the final time one: i.e., the probability that $h_n \geq 1 - \e$
is small when $\e > 0$ is small, uniformly in $n$.
Suppose now that late coalescence is indeed absent, and consider the random difference $\wgt{n}{x}{z+ \eta} - \wgt{n}{x}{z}$ 
in the cases that $x = -1$ and $x = 1$. As Figure~\ref{f.heuristics}(left) illustrates,
when $\vert \eta \vert$ small enough, this difference may be expected to be shared between the two cases, forcing the weight difference profile to be locally constant near $z$. Indeed, as  
Figure~$1$(a) illustrates, the difference in polymer trajectories between $(-1,0) \to (z + \eta,1)$
and $(-1,0) \to (z,1)$ is given by a diversion of trajectory only after the coalescence time $h_n(z)$; and this same polymer trajectory difference holds  between $(1,0) \to (z + \eta,1)$
and $(1,0) \to (z,1)$.

\subsubsection{Heuristics $3$: the Hausdorff dimension of exceptional points.}
Proving the lower bound on the Hausdorff dimension of a random fractal is often more demanding than deriving the upper bound. In this problem, however, two seemingly divorced considerations yield matching lower and upper bounds. Local Gaussianity of weight profiles forces the dimension to be at least one-half; while the rarity of disjoint pairs of polymers with nearby endpoints yields the matching upper bound.

\begin{figure}
\includegraphics[width=0.6\textwidth]{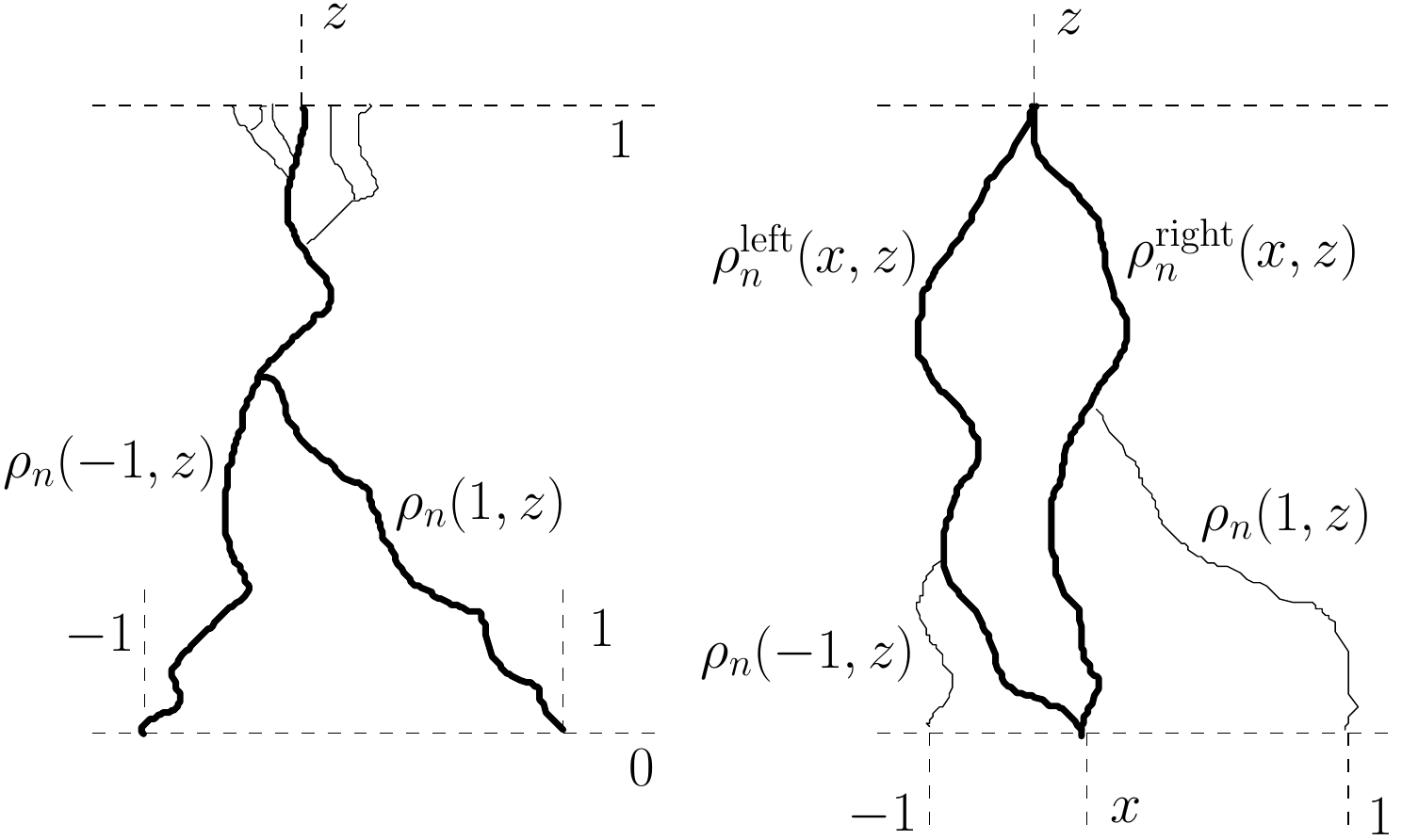}
\caption{{\em Left.} When $z \in \R$ is given, the bold curves $\rho_n(-1,z)$ and $\rho_n(1,z)$
coalesce at a random height which has a probability of being close to one that is low, uniformly in high $n$. The trunks of the trees rooted at $-1$ and $1$ being shared above an intermediate height, the trees' canopies are also shared around the tip~$(z,1)$ of the trunk; heuristically at least follows the local constancy of $y \to W_n(1,y) - W_n(-1,y)$ near $z$.  
  {\em Right.} In this case, $z \in E'$; the bold curves are $\rho_n^{\rm left}(x,z)$ and $\rho_n^{\rm right}(x,z)$; and the fainter curves, namely $\rho_n(\pm 1,z)$, merge with the respective bold curves as height rises. Since $z \in E'$, the curves $\rho_n(\pm 1,z)$ are special in that they meet only at height one; some $x \in [-1,1]$ exists for which the polymers $\rho_n^{\rm left}(x,z)$ and $\rho_n^{\rm right}(x,z)$ have the doubly remarkable characteristic that intersection occurs merely at heights zero and one despite both endpoints being shared.}\label{f.heuristics} 
\end{figure}

There are thus two tasks that require overview. 
For the lower bound, take $\e > 0$ small and let $\mathscr{C}_\e$ be the set of $\e$-length subintervals of $[-1,1]$ of the form $[\e k, \e (k+1)]$, $k \in \Z$, on which $Z$ is {\em not} constant. Since $Z$ is a difference of Airy$_2$ processes -- formally, $Z(v)$ equals $\wgt{n}{-1}{v} - \wgt{n}{1}{v}$  with $n = \infty$ -- and the Gaussian local variation of these processes is gauged by Theorem~\ref{t.differenceweight}, $Z$ may vary on a length-$\e$ interval only by order $\e^{1/2}$. Since  $Z(1) - Z(-1)$ is a random quantity of unit order, we see that typically $\vert \mathscr{C}_\e \vert \geq O(1) \e^{-1/2}$, whence, roughly speaking, is the Hausdorff dimension of $E$ seen to be at least one-half.

Deriving the matching upper bound is a matter of showing that the exceptional set of points $E$ about whose members $Z$ is not locally constant is suitably sparse. In view of the preceding argument for the almost everywhere local constancy of $Z$, we see that $E \subseteq E'$,  
where $E'$ denotes the set of  $z \in \R$ such that the paths $\rho_n(-1,z)$ and $\rho_n(1,z)$ coalesce only at the final moment, at height one (when we take $n = \infty$ formally). The plan is to argue that the number $\e$-length intervals in a mesh that contain such a point $z$ is typically of order at most $\e^{-1/2}$, since then an upper bound on the Hausdorff dimension of $E \subseteq E'$ follows directly. Suppose that $z \in E'$
and consider dragging the lower spatial endpoint~$u$ of the polymer $\polymer{n}{u}{z}$
rightwards from $u = -1$ until the first moment~$x$ at which the moving polymer intersects its initial condition $\polymer{n}{-1}{z}$ only at the ending height~one -- see Figure~\ref{f.heuristics}(right) for a depiction. That $z \in E'$ implies that $x \leq 1$. The journey $(x,0) \to (z,1)$ is doubly special, since polymer disjointness is achieved at both the start and the end of the journey. 
Indeed, there exists a pair of polymers, which may be called
$\rho_n^{{\rm left}}(x,z)$
 and 
$\rho_n^{{\rm right}}(x,z)$
each running from $(x,0)$ to $(z,1)$, that are disjoint except at these endpoints.
Consider intervals $I$ and $J$ between consecutive elements of the mesh $\e \Z$
that respectively contain $x$ and $y$. 
The event $\maxpolymac{n}{I}{J} \geq 2$ occurs, because the just recorded polymer pair {\em in essence} realizes it; merely in essence due to the meeting at start and end, a problem easily fixed. 
Theorem~\ref{t.disjtpoly.pop} shows that the dominant order of this event's probability is at most~$\e^{3/2}$.
 In light of this bound, the total number of such interval-pairs $I \times J$ inside $[-1,1] \times [-1,1]$ is at most $\e^{-1} \cdot \e^{-1} \cdot \e^{3/2} = \e^{-1/2}$. 
Since each $\e$-interval in a mesh that intersects $E'$ furnishes a distinct such pair $(I,J)$, we see that such intervals typically number at most order $\e^{-1/2}$, as we sought to show.

\subsection{Sharpness of the estimate on rarity of polymer pairs with nearby endpoints}

Conjecture~$1.3$ in~\cite{brownianLPPtransversal} asserts that Theorem~\ref{t.disjtpoly.pop} is sharp in the sense that no improvement can be made in the exponent $k^2 -1$. The method of proof of Theorem~\ref{t.z}
proves 
this conjecture in the case that $k=2$.  

\begin{theorem}\label{t.disjtpoly.lb} 
There exists $d > 0$ such that, for $\eta > 0$, we may find $\e_0 = \e_0(d,\eta)$
for which, whenever  $\e \in (0,\e_0)$, there exists 
$n_0 = n_0(d,\e,\eta)$ so that $n \in \N$, $n \geq n_0$, implies that
$$
     \PP \Big( \maxpolymac{n}{[0,2\e]}{[0,2\e]}  \geq 2 \Big) 
  \, \geq  \,
 d \, \e^{3/2 + \eta} \, .
$$
\end{theorem}
This result implies directly that 
$$
 \limsup_{\e \searrow 0} \, \limsup_n \, \frac{\log \PP \Big( \maxpoly_n \big( [0,2\e] , [0,2\e] \big)  \geq k \Big)}{\log \e}  \, \leq \, \frac{k^2 - 1}{2} 
 $$
 when $k = 2$; after the double replacement of $[0,2\e]$ by $[-\e,\e]$ -- replacements permitted by the stationary increments of the underlying noise field $B$ -- we indeed obtain~\cite[Conjecture~$1.3$]{brownianLPPtransversal} with~$k~=2$.

\subsection{Polymer basics}

A splitting operation on polymers will be needed.
\begin{definition}\label{d.split}
Let $n \in \N$ and let $x,y \in \R$ verify $y - x \geq 2^{-1} n^{1/3}$. Let $\rho$ denote a polymer from $(x,0)$ to $(y,1)$, and let $(z,s) \in \R \times [0,1]$ be an element of $\rho$ for which $s \in n^{-1}\Z$; in this way, $z$ lies in one of $\rho$'s horizontal planar line segments. The removal of $(z,s)$ from $\rho$ generates two connected components. Taking the closure of either of these amounts to adding the point $(z,s)$ back to the component in question. The resulting sets are $n$-zigzags from $(x,0)$ to $(z,s)$ and from $(z,s)$ to $(y,1)$, and it is a simple matter to check that each of these zigzags is in fact a polymer. Denoting these two polymers by $\rho_-$ and $\rho_+$, we use the symbol $\circ$ evoking concatenation to express this splitting of $\rho$ at $(z,s)$, writing $\rho = \rho_- \circ \rho_+$.
\end{definition}
 
We have mentioned that \cite[Lemma~$4.6(1)$]{hammond2017patchwork} implies that the polymer making the journey $(x,0)$ to $(y,1)$ is almost surely unique for any given $x,y \in \R$ for which it exists; namely, for those $(x,y)$ satisfying $y - x \geq - 2^{-1} n^{1/3}$. Although it may at times aid intuition  to consider the almost surely unique such polymer  $\rho_n(x,y)$, as we did in the preceding heuristical presentation, it is not logically necessary for the presentation of our proofs, which we have formulated without recourse to almost sure polymer uniqueness. As a matter of convenience, we will sometimes invoke the almost sure {\em existence} of polymers with given endpoints; this result is an exercise that uses compactness and invokes the continuity of the underlying Brownian ensemble $B:\Z \times \R \to \R$.  

A few very straightforward properties of zigzags and polymers will be invoked implicitly: examples include that any pair of zigzags that intersect do so at a point, necessarily of the form $(u,s) \in \R \times n^{-1}\Z$, that lies in a horizontal line segment of both zigzags; and that the subpath of a polymer between two of its members having this form is itself a polymer. 

\section{The proofs of the main theorems}

By far the hardest element of Theorem~\ref{t.z} is its third assertion, concerning Hausdorff dimension. After introducing a little notation and recalling the definition of this dimension, we reformulate Theorem~\ref{t.z}(3) as the two-part Theorem~\ref{t.hd} in which the needed upper and lower bounds are expressed. These bounds are then proved in ensuing two subsections.
A fourth subsection provides the proof of Theorem~\ref{t.disjtpoly.lb}. 

We set $Z_n: \R \to \R$ to be the weight difference profile 
$$
 Z_n(z) = \wgt{n}{1}{z} - \wgt{n}{-1}{z} \, , 
$$
where the domain of definition of $Z_n$ may be chosen to be $\R$ by use of the convention specified before Theorem~\ref{t.z}. Recall from the theorem and for use shortly that $Z$ denotes any weak limit point of the random functions $Z_n$.

Let $f$ be a real-valued function defined on $\R$ or a compact interval thereof. We will write $\mathrm{LV}(f)$ for the subset of the domain of $f$ that comprises points $z$ of {\em local variation} of $f$ about which no interval exists on which $f$ is constant.

\begin{definition}\label{d.hausdorff}
Let $d\in [0,\infty)$. 
The $d$-dimensional Hausdorff measure $H^d(X)$ of a metric space~$X$ equals $\lim_{\delta\searrow 0} H^d_\delta(X)$ where, for $\delta>0$, we set 
$$
H^{d}_{\delta}(X) \, = \, \inf \, \bigg\{  \, \sum_{i} \mathrm{diam}(U_i)^d: \{U_i\}~\text{is a countable cover of}~X~\text{with}~0< \mathrm{diam}~U_{i}<\delta \, \bigg\} \, .
$$
The Hausdorff content $H^d_\infty(X)$ of $X$ is specified to by choosing $\delta = \infty$ here, a choice that renders vacuous the diameter condition on the covers.


The Hausdorff dimension $d_H(X)$ of $X$ equals the infimum of those positive $d$ for which $H^d(X)$ equals zero; and it is straightforwardly seen that the Hausdorff measure $H^d(X)$ may here be replaced by the Hausdorff content $H^d_\infty(X)$
to obtain an equivalent definition.
\end{definition}

We will write $\vert U_i \vert$ in place of $\mathrm{diam}~U_i$, doing so without generating the potential for confusion because every considered $U_i$ will be an interval.

{\bf Proof of Theorem~\ref{t.z}: (1).} 
By \cite[Lemma 2.2(1)]{hammond2017modulus}, for each $n \in \N$ and $x \in \R$, the random function $z \to \wgt{n}{x}{z}$
is almost surely continuous on its domain of definition $z \geq x - 2^{-1} n^{1/3}$.
The process $Z: \R \to \R$ is thus a weak limit point of continuous stochastic processes mapping the real line to itself. The Skorokhod representation of weak convergence thus implies that the prelimiting processes may be coupled with the limit $Z$ in such a way that, almost surely, they converge locally uniformly to $Z$. Thus $Z$ is seen to be continuous almost surely.

To show that $Z:\R\to \R$ is non-decreasing, we will derive a counterpart monotonicity assertion in the prelimit. Indeed, it will be enough to argue that
$$
 Z_n(z) \geq Z_n(y) \, \, \, \textrm{whenever $y,z \in \R$, $y < z$ and $n \in \N$ satisfy $2^{-1} n^{1/3} \geq \max\{\vert y \vert, \vert z \vert \} + 1$} \, .
$$
The indicated inequality on parameters is needed merely to ensure that the concerned weights $W_n(\pm 1, v)$, $v \in \{ y,z\}$, are well specified by the defining formula~(\ref{e.weight})
With the condition imposed, there almost surely exist polymers, which we denote by $\rho^1$ and $\rho^2$, that make the respective journeys $(-1,0) \to (z,1)$ and $(1,0) \to (y,1)$.
 By planarity, we may find an element $(w,s) \in \R \times [0,1]$ of $\rho^1 \cap \rho^2$ with $s \in n^{-1}\Z$. Let $\rho^1 = \rho^1_- \circ \rho^1_+$ and $\rho^2 = \rho^2_- \circ \rho^2_+$ denote the polymer decompositions resulting from splitting the two polymers at $(w,s)$. We write $W^i_{\pm}$ with $i \in \{1,2\}$ for the weights of the four polymers so denoted.

Note that $W_n(-1,z) = W^1_- + W^1_+$ and  $W_n(1,y) = W^2_- + W^2_+$. The quantity $W_n(1,z)$ is at least the weight of $\rho^2_- \circ \rho^1_+$; which is to say, $W_n(1,z) \geq W^2_-  + W^1_+$.
Likewise, $W_n(-1,y) \geq W^1_- + W^2_+$. We have two equalities and two inequalities -- we use them all to prove the bound that we seek.
Indeed, we have that 
$$Z_n(z) = W_n(1,z) - W_n(-1,z) \geq \big( W^2_-  + W^1_+ \big) - \big( W^1_- + W^1_+ \big) = W^2_- - W^1_-.$$ We also see that $$Z_n(y) = W_n(1,y) - W_n(-1,y) \leq \big( W^2_- + W^2_+ \big) - \big( W^1_- + W^2_+ \big) = W^2_- - W^1_-.$$ That is, $Z_n(z) \geq Z_n(y)$, as we sought to show. 

{\bf (2).} This is implied by the third part of the theorem.

{\bf (3).} This follows from the next theorem. \qed

\begin{theorem}\label{t.hd}
\begin{enumerate}
\item The Hausdorff dimension of $\mathrm{LV}(Z)$ is at most one-half almost surely.
\item 
Let $\delta>0$. There exists $M=M(\delta)>0$ such that, with probability at least $1-\delta$, $\mathrm{LV}(Z)\cap [-M,M]$ has Hausdorff dimension at least one-half. 
\end{enumerate}
\end{theorem}

\subsection{The upper bound on Hausdorff dimension}

Here we prove Theorem~\ref{t.hd}(1). 
The principal component is the next result, which offers  control on  the $d$-dimensional Hausdorff measure of $\mathrm{LV}(Z_n)$ for $n$ finite but large.
The result is stated for the prelimiting random functions $Z_n$
in order to quantify explicitly the outcome of our method, but, for our application, we want to study the weak limit point $Z$. With this aim in mind, we present Theorem~\ref{t.ub}, and further results {\em en route} to Theorem~\ref{t.hd}(1), so that assertions are made about both the prelimit and the limit. The notational device that permits this to set $Z_\infty$ equal to the weak limit point $Z$; thus the choice $n=\infty$ corresponds to the limiting case.

\begin{theorem}\label{t.ub}
Let $d> 1/2$ and $M>0$. Consider any positive sequences $\big\{ \delta_k: k \in \N \big\}$ and  $\big\{ \eta_k: k \in \N  \big\}$  that converge to zero.  For each $k \in \N$, there exists $n_k=n_k \big( d,M,\delta_{k}, \eta_{k} \big)$ such that, for $n \in \N \cup \{ \infty\}$ with $n \geq n_k$,
$$
\PP \Big( H^{d}_{\infty}\big( \mathrm{LV}(Z_n) \cap [-M,M] \big) \leq \eta_{k} \Big) \geq 1 - \delta_k \, .
$$
\end{theorem}

\begin{lemma}\label{l.localconstancy}
Let $n \in \N$ and let $z_1,z_2 \in \R$ satisfy $z_1 < z_2$. Suppose that there exist  polymers
making the journeys $(-1,0) \to (z_1,1)$ and $(1,0) \to (z_2,1)$
whose intersection is non-empty. Then $Z_n$ is constant on $[z_1,z_2]$. 
\end{lemma}
{\bf Proof.} For two zigzags $\zeta$ and $\zeta'$ that begin and end at respective times zero and one, we write $\zeta \preceq \zeta'$ to indicate that `$\zeta'$ is on or to the right of $\zeta$', in the sense that $\zeta'$ is contained in the union of the semi-infinite horizontal planar line segments whose left endpoints are elements of $\zeta$.

Let $\rho_{\rm left}$ and $\rho_{\rm right}$ be polymers of respective journeys $(-1,0) \to (z_1,1)$ and $(1,0) \to (z_2,1)$ whose existence is hypothesised. 
Let $(u,s) \in \R \times [0,1]$ with $s \in n^{-1} \Z$ denote an element of  $\rho_{\rm left} \cap \rho_{\rm right}$. 
Our {\em first claim} is that we may impose that $\rho_{\rm left} \preceq \rho_{\rm right}$
while respecting all of these properties. To verify this, note that, should this ordering condition fail, 
$\rho_{\rm right}$ makes at least one excursion to the left of $\rho_{\rm left}$, in the sense that there exist a pair of elements in these two polymers whose removal from each results in a pair of zigzags that connect the pair, with the one arising from $\rho_{\rm left}$ lying on or to the right of that arising from $\rho_{\rm right}$. The weight of these two zigzags is equal, and each may be recombined with the remaining subpaths of the opposing polymer to form updated copies of  
$\rho_{\rm left}$ and $\rho_{\rm right}$ in which the excursion in question has been eliminated. There are only finitely many excursions, because the vertical intervals assumed by excursions are disjoint and abut elements of $n^{-1}\Z \cap [0,1]$. Thus, after finitely many iterations of this procedure, will the condition  $\rho_{\rm left} \preceq \rho_{\rm right}$ be secured. Any point that changes hands in the course of this operation does so because it belongs to exactly one of the original copies of  $\rho_{\rm left}$ and $\rho_{\rm right}$. Since $(u,s)$ lies in the intersection of them,  it remains in the intersection at the end of the procedure. Thus is this first claim verified.

Let $z \in [z_1,z_2]$, and let $\rho$ denote a polymer from $(-1,0)$ to $(z,1)$. By a {\em second claim}, we may impose the sandwiching condition that $\rho_{\rm left} \preceq \rho \preceq  \rho_{{\rm right}}$.
Indeed, any excursion that $\rho$ makes to the left of $\rho_{{\rm left}}$ may be substituted by the intervening trajectory of that polymer; and likewise for $\rho_{{\rm right}}$; so that this second claim is seen to hold.

Consistently with the use of notation~$\preceq$, two closed horizontal planar intervals  $A$ and $B$ at a given height verify $A \preceq B$
when the respective endpoints to $A$ are at to the left of those of $B$. Indeed, we have that 
$$
\rho_{\rm left} \cap \big(  \R \times \{ s \} \big) \, \preceq \, \rho \cap \big(  \R \times \{ s \} \big) \, \preceq \,  \rho_{{\rm right}} \cap \big(  \R \times \{ s \} \big) \, .
$$  
Since the first and third horizontal planar intervals contain $(u,s)$, we see that $\rho$ also contains $(u,s)$.

We consider the decomposition  $\rho = \rho_- \circ \rho_+$ from Definition~\ref{d.split}, where $\rho$ is split at $(u,s)$. Similarly we denote $\rho_{\rm left} = \rho_{{\rm left},-} \circ \rho_{{\rm left},+}$
and $\rho_{{\rm right}} = \rho_{{\rm right},-}  \circ \rho_{{\rm right},+}$, with the splits again occurring at~$(u,s)$.

Our {\em third claim} is that $\rho_{{\rm left},-}  \circ \rho_+$ is a polymer from $(-1,0)$ to $(z,1)$; and that $\rho_{{\rm right},-} \circ \rho_+$ is a polymer from $(1,0)$ to $(z,1)$. Indeed, 
the weight of $ \rho_{{\rm left},-}  \circ \rho_+$ is the sum of the weights of its constituent paths, of which the first is at least the weight of $\rho_-$, since $\rho_{{\rm left},-}$ is a polymer that shares its endsponits with $\rho_-$.  Thus the weight of $\rho_{{\rm left},-}  \circ \rho_+$  is at least that of $\rho$. Since $\rho$ is a polymer that shares the endpoints of  $\rho_{{\rm left},-}  \circ \rho_+$, the latter zigzag is a polymer.  The second element in the third claim is similarly verified.

The quantity $Z_n(z)$, being $\wgt{n}{1}{z} - \wgt{n}{-1}{z}$, equals
$W_1 - W_2$, where $W_1$ 
 the sum of the weights of $\rho_{{\rm right},-}$ and $\rho_+$; and $W_2$ is the sum of the weights of $\rho_{{\rm left},-}$
 and $\rho_+$. Since $W_1 - W_2$, being the difference in weight between $\rho_{{\rm right},-}$ and $\rho_{{\rm left},-}$, is independent of $z \in [z_1,z_2]$, the proof of Lemma~\ref{l.localconstancy} is complete. \qed

\begin{proposition}\label{p.drag}
Let $n \in \N$, $z \in \R$ and $\e > 0$. 
Suppose that there is no point of intersection between any pair of polymers making the respective journeys $(-1,0) \to (z,1)$ and $(1,0) \to (z+\e,1)$. 
Then there exists an interval $I \subset [-1,1]$
of length~$\e$ for which $\maxpolymac{n}{I}{[z,z+\e]} \geq 2$.
\end{proposition}
{\bf Proof.} 
For $x_1, x_2, z_1, z_2 \in \R$ for which $x_1 < x_2$ and $z_1 < z_2$,   let $\mathrm{NonInt}_{n}\big( \{x_1,x_2\} , \{z_1,z_2\} \big)$ denote the event that  there is no point of intersection between any polymer from $(x_1,0)$ to $(z_1,1)$ and any polymer from $(x_2,0)$ to $(z_2,1)$.

Let $X$ denote the supremum of those $x \in [-1,1]$ for which $\mathrm{NonInt}_n \big( \{x,1\} , \{z,z + \e\} \big\}$ occurs. 
 Since $x=-1$ qualifies, $X$ is a well-defined random variable, taking values in $[-1,1]$.
  We will first treat the trivial case that $X = 1$ and then turn to the principal one, when  $X \in [-1,1)$.

When $X = 1$, we may find $u \in (1 - \e,1)$
for which  $\mathrm{NonInt}_n \big( \{u,1\} , \{z,z + \e\} \big)$ occurs. 
Thus Proposition~\ref{p.drag} holds with $I = [1 - \e, 1]$.

Suppose instead then that $X \in [-1,1)$.  If we further insist that $X > -1$, we may locate $u > -1$, $u \in (X - \e,X)$ and $v < 1$, $v < u+\e$, for which the event $\mathrm{NonInt}_n \big( \{x,1\} , \{z,z + \e \} \big)$ occurs when $x = u$ and does not occur when $x = v$. If, on the other hand, $X$ is equal to $-1$, 
we may achieve the same circumstance by taking $u = -1$ and $v = u+\e$.

Consider any polymer $\rho_1$ from $(u,0)$ to $(z,1)$, and note that $\rho_1$ is disjoint from any polymer from $(1,0)$ to $(z + \e,1)$. 
We may, by definition, find polymers  $\rho_2$ and $\rho_3$ of non-empty intersection that make the respective journeys $(v,0) \to (z,1)$ and $(1,0) \to (z + \e,1)$. Let  $(w,t) \in \R \times [0,1]$, with $t \in n^{-1}\Z$, be an element of $\rho_2 \cap \rho_3$.

Write $\rho_2 = \rho_{2,-} \circ \rho_{2,+}$ where the right-hand polymers are formed by splitting $\rho_2$ at $(w,t)$; and use the counterpart notation $\rho_3 = \rho_{3,-} \circ \rho_{3,+}$.

\begin{figure}[h]
\centering
\includegraphics[scale=.7]{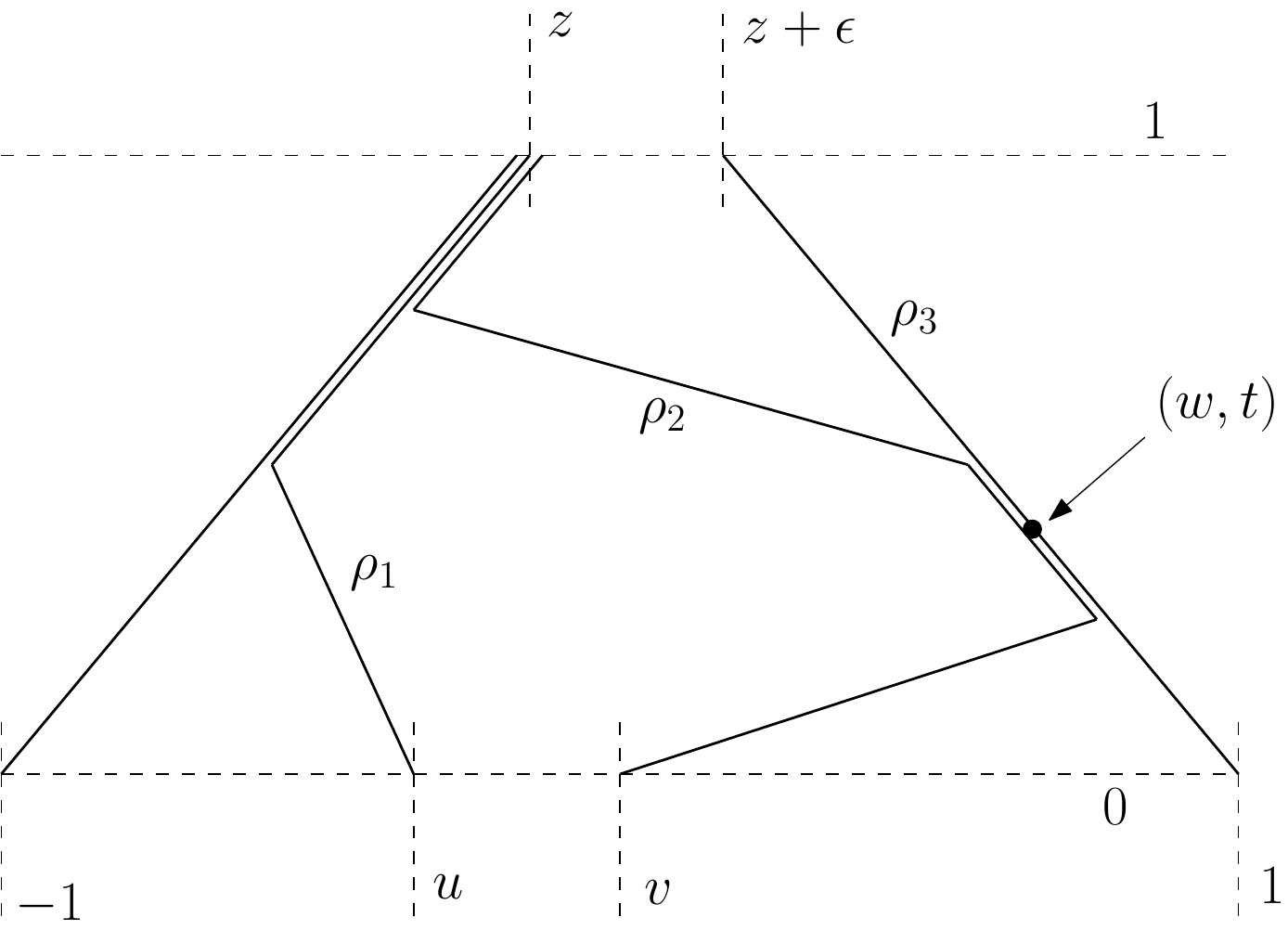}
\caption{Four polymers are depicted in a caricature as piecewise affine curves; close parallel line segments indicate shared sections among these curves. The leftmost, unlabelled, curve is a polymer from $(-1,0)$ to $(z,1)$. The path $\rho$, which we claim to be a polymer, follows $\rho_2$ from $(v,0)$ to $(w,t)$, from where it follows $\rho_3$ to its end at $(z+\e,1)$.}
\label{f.fourcurvecaricature}
\end{figure}

Consider the path $\rho = \rho_{2,-} \circ \rho_{3,+}$. We {\em claim} that $\rho$ is a polymer from $(v,0)$ to $(z+\e,1)$ -- see Figure~\ref{f.fourcurvecaricature}.

We will establish this by considering any polymer $\rho'$
 from $(v,0)$ to $(z+\e,1)$, and arguing that the weight of $\rho$ is at least that of $\rho'$. 
 
 This claim will be proved via an intermediate step, in which we exhibit a polymer $\bar\rho$ from $(v,0)$ to $(z+\e,1)$ such that  
 $(w,t) \in \bar\rho$. To construct $\bar\rho$, we indeed consider any polymer $\rho'$
 from $(v,0)$ to $(z+\e,1)$. By planarity, $\rho'$ intersects either $\rho_{3,-}$
 or $\rho_{2,+}$. Suppose that $\rho' \cap \rho_{3,-} \not= \emptyset$. Set $\bar\rho$ equal to the zigzag formed by following $\rho'$ until its first intersection with $\rho_{3,-}$, and then following $\rho_3$ until its end. Then $\bar\rho$ runs from $(v,0)$ to $(z+\e,1)$; $(w,t) \in \bar\rho$; and, since the weight of that part of $\bar\rho$ that runs along $\rho_3$ is at least the weight of that part of $\rho'$ that runs over the journey with the same endpoints in light of $\rho_3$ being a polymer, we see that $\bar\rho$ is itself a polymer. If instead  $\rho' \cap \rho_{2,+} \not= \emptyset$, then a suitable~$\bar\rho$ may be formed by running along $\rho_{2,-} \circ \rho_{2,+}$ until an element of  $\rho' \cap \rho_{2,+}$ is encountered, and then following the course of $\rho'$ to its end at $(z+\e,1)$.
 
The zigzag $\rho$ makes its journey in two stages, pausing at $(w,t)$ between them. Likewise for the polymer $\bar\rho$. But each stage for $\rho$ is a polymer, so the weight of each stage for $\rho$ must be at least what it is for $\bar\rho$. Thus we see that the weight of $\rho$ is at least that of $\bar\rho$, so that we confirm the claim that $\rho$ is a polymer.

We now argue that $\rho$ and $\rho_1$ are disjoint.
With a view to obtaining a contradiction, suppose instead that $\rho_1 \cap \rho \not= \emptyset$. Since $\rho = \rho_{2,-} \circ \rho_{3,+}$, a point of intersection must lie in either $\rho_{2,-}$ or $\rho_{3,+}$. The latter is impossible, because in that case, we would find a point in $\rho_1 \cap \rho_3$, and none exists since  $\mathrm{NonInt}_n \big( \{u,1\} , \{z,z + \e \} \big)$ occurs. On the other hand, were $\rho_1 \cap \rho_{2,-}$ non-empty, then we might take a journey along $\rho_1$ until its first intersection with $\rho_{2,-}$, and then follow the course of~$\rho_2$ until its end at $(z,1)$. Since the part of $\rho_2$ that is followed is a polymer, this journey is itself a polymer from $(u,0)$ to $(z,1)$. But the journey visits $(w,t) \in \rho_3$, in violation of the occurrence of   $\mathrm{NonInt}_n \big( \{u,1\} , \{z,z + \e \} \big)$.

Since    $\rho \cap \rho_1 = \emptyset$, we see that  $\big\{ \maxpolymac{n}{[u,v]}{[z,z+\e]} \geq 2 \big\}$ occurs. By setting~$I$ equal to any interval of length~$\e$ such that $[u,v] \subseteq I \subseteq [-1,1]$, Proposition~\ref{p.drag} has been obtained in the case that $X \in [-1,1)$.  \qed

\begin{proposition}\label{p.mdp}
Let $n \in \N$, $z \in \R$ and $\e > 0$.
When the event that $\mathrm{LV}(Z_n) \cap [z,z + \e] \neq \emptyset$ occurs, there almost surely exists  $u \in [-1,1-2\e] \cap \e\Z$ such that  $\maxpolymac{n}{[u,u+2\e]}{[z,z+\e]} \geq 2$. 
\end{proposition}
{\bf Proof.}
 By Lemma~\ref{l.localconstancy} and Proposition~\ref{p.drag}, the occurrence of $\mathrm{LV}(Z_n) \cap [z,z + \e] \neq \emptyset$ entails the existence of $u\in [-1,1-\e]$
for which $\maxpolymac{n}{[u,u+\e]}{[z,z+\e]} \geq 2$. If we replace the interval $[u,u+\e]$ by $[u,u+2\e]$ (and $u$ by $u - \e$ should $u$ be at least $1-2\e$), we may further demand, as we need to do in order to prove the proposition, that $u \in [-1,1-2\e] \cap \e\Z$. \qed

\begin{proposition}\label{p.localnonconst}
There exists a positive constants $C_0$ and $C_1$ such that, for $M > 0$, we may find  $\e_0 = \e_0(M) > 0$  for which, when $\e \in (0,\e_0)$ and $n \in \N \cup \{ \infty \}$ satisfies $n \geq C_1 (M+2)^{36}\e^{-C_1}$, we have that
$$ 
\PP \Big( \mathrm{LV}(Z_n) \cap [z,z + \e] \neq \emptyset \Big) \, \leq \,  \e^{1/2} \cdot \exp \Big\{ C_0 \big( \log \e^{-1} \big)^{5/6} \Big\} 
$$ 
whenever $z \in [-M,M-\e]$.
\end{proposition}
{\bf Proof.} We defer consideration of the case that $n = \infty$ and suppose that $n \in \N$.
By Proposition~\ref{p.mdp}, the event whose probability we seek to bound above is seen to entail the existence of a pair of disjoint polymers that make the journey $[u,u+2\e] \to [z,z+\e]$ between times zero and one. 
Theorem~\ref{t.disjtpoly.pop} with $k=2$, $x = u + \e$ and $y = z$ 
provides an upper bound on the probability of this polymer pair's existence for given~$u$, since the condition that $n \geq C_1 (M+2)^{36}\e^{-C_1}$ for a suitably high choice of the constant $C_1$ permits the use of this theorem.
A union bound over the at most $2\e^{-1}$ choices of $u$ provided by the use of Proposition~\ref{p.mdp} then
yields Proposition~\ref{p.localnonconst} for finite choices of~$n$, where suitable choices of $\e_0$ and $C_0$ absorb the factor of $2\e^{-1}$ generated by use of the union bound.

To treat the case that $n = \infty$, note that, by the Skorokhod representation, the processes $Z_n$ indexed by finite $n$ may be coupled to the limit $Z_\infty$ so as to converge along a suitable subsequence uniformly on any compact set. Momentarily relabelling so that $Z_n$ denotes the convergent subsequence, it follows that, for any closed interval $I \subseteq [-M,M]$, 
$$
\limsup_{n\to \infty}\P\big( Z_n~\text{is constant on}~I \big) \, \leq \, \P \big( Z~\text{constant on}~I \big) \, .
$$ 
Thus does
Proposition \ref{p.localnonconst} 
in the remaining case that $n = \infty$ follow from the case of finite~$n$. \qed

For given $M > 0$, let $N_n(M)$ denote the number of intervals $[u,u+\e]$ with $u \in \e \Z$ that intersect $[-M,M]$
and on which $Z_n$ fails to be constant. Proposition~\ref{p.localnonconst} permits us to bound the upper tail of $N_n(M)$.

\begin{corollary}\label{c.number}
There exists $n_0 = n_0(\e,M)$ such that, for $n \in \N \cup \{ \infty \}$ with $n \geq n_0$,
$$ 
\PP \left( N_n(M) \geq\e^{-1/2} \cdot   4\delta^{-1}  M  \exp \Big\{ C_0 \big(\log \e^{-1} \big)^{5/6} \Big\} \right) \, \leq \, \delta \, .
$$
\end{corollary}
{\bf Proof.}
Proposition~\ref{p.localnonconst} implies that $\E N_n(M) \leq 4M \e^{-1/2} \exp \big\{ C_0 (\log \e^{-1})^{5/6} \big\}$, so that Markov's inequality implies the desired result. \qed

{\bf Proof of Theorem \ref{t.ub}.}
Recall that $d> 1/2$ and $M > 0$; and that $\big\{ \delta_k: k \in \N \big\}$ and $\big\{ \eta_k: k \in \N \big\}$ are positive sequences that are arbitrary subject to their converging to zero. Let $k \in \N$. We must, on an event of probability at least $1-\delta_k$, exhibit for all $n \in \N \cup \{ \infty \}$ verifying $n \geq n_k$ a cover of $\mathrm{LV}(Z_{n})\cap [-M,M]$ comprised of intervals $U_{i}$ that satisfy $\sum_{i} |U_i|^{d}\leq \eta_k$. Here, $n_k$ may depend on $d$, $M$, $\delta_k$  and $\eta_k$. 

The cover is chosen to be equal to the set of those intervals of the form $[u,u+\e]$ with $u \in \e \Z$
whose intersection with  $\mathrm{LV}(Z_{n}) \cap [-M,M]$ is non-empty.
Corollary~\ref{c.number} implies that it is with probability at least $1-\delta_k$ that
$$
\sum_i |U_i|^d \leq 4 \delta_{k}^{-1} M C_0 \e^{d-1/2} \exp \big\{ C_0 (\log \e^{-1})^{5/6} \big\} \, ,
$$
provided that $n$ exceeds a value that is determined by $M$ and $\e$. Since $d> 1/2$, this right-hand side converges to zero in the limit of $\e \searrow 0$ provided that every other parameter is held fixed. Recalling the given sequences $\delta$ and $\eta$, we may select $\e_0 = \e_0\big(M,d,\delta_k,\eta_k\big)$ so that, when $\e \in (0,\e_0)$, the preceding right-hand side is at most $\eta_k$ whenever the parameter $n$ to chosen to be high enough. Thus do we conclude the proof of Theorem~\ref{t.ub}. \qed

{\bf Proof of Theorem \ref{t.hd}(1).}
It follows directly from Theorem \ref{t.ub} with $n = \infty$ that, given any $d > 1/2$; any summable sequence $\big\{ \delta_k: k \in \N \big\}$; any sequence $\big\{ \eta_k : k \in \N \big\}$ that converges to zero; and further any sequence $\big\{M_{k}: k \in \N \big\}$ that converges to $\infty$;  there exists,  with probability at least $1 - \delta_k$,  a countable cover of $\mathrm{LV}(Z_\infty)\cap [-M_k,M_k]$ that witnesses the Hausdorff content $H^d_\infty \big( \mathrm{LV}(Z_\infty)\cap [-M_k,M_k] \big)$ being less than $\eta_k$. A use of the Borel-Cantelli lemma then shows that almost surely there exists a random $K_0 \in \N$ such that, for $k \geq K_0$, $H^d_\infty\big(\mathrm{LV}(Z_\infty)\cap [-M_k,M_k] \big) \leq \eta_k$. Since $\eta_k$ converges to zero, we see that $H^d_\infty \big( \mathrm{LV}(Z_\infty) \big)$ is zero almost surely for every $d> 1/2$; and thus do we prove Theorem~\ref{t.hd}(1).  \qed

\subsection{The matching lower bound}

Here we prove Theorem~\ref{t.hd}(2).
We will do so by invoking the following {\em mass distribution principle}, a tool that offers a lower bound on the Hausdorff dimension of a set which supports a non-trivial measure that attaches low values to small balls. In this assertion, a mass distribution is a measure $\mu$ defined on the Borel sets of a metric space $E$ for which $\mu(E) \in (0,\infty)$.

\begin{theorem}\cite[Theorem~$4.19$]{PeresMorters}\label{t.mdp}
Suppose given a metric space $E$ and a value $\alpha > 0$. For any mass distribution $\mu$ on $E$, and any positive constants $K$ and $\eta > 0$, the condition that 
\begin{equation}\label{e.vcdelta}
\mu(V) \leq K \vert V \vert^\alpha 
\end{equation}
for all closed sets $V \subseteq E$ of diameter $\vert V \vert$ at most $\eta$
ensures that the Hausdorff measure $H^\alpha(E)$ is at least $K^{-1} \mu(E) > 0$; and thus that the Hausdorff dimension $d_H(E)$ is at least $\alpha$.
\end{theorem}
The set $\mathrm{LV}(Z)$ under study in Theorem~\ref{t.hd}(2) supports a natural random measure $\mu$ in view of Theorem~\ref{t.z}(1): we may specify  $\mu(a,b] = Z(b) - Z(a)$ for $a,b \in \R$ with $a \leq b$, so that $Z$ is the distribution function of $\mu$.

For $M > 0$, we aim to apply Theorem~\ref{t.mdp} for any given $\alpha  \in (0,1/2)$, with $E = [-M,M]$ and $\mu$ given by restriction to $E$. What is needed are two inputs: an assertion of  {\em non-degeneracy} that the so defined $\mu$ is typically positive 
when $M > 0$ is high; and an assertion of {\em distribution of measure} --  absence of local concentration for $\mu$ -- that will validate the hypothesis~(\ref{e.vcdelta}). 
  
  We present these two inputs; use them to prove Theorem~\ref{t.hd}(2) via Theorem~\ref{t.mdp}; and then prove the two input assertions in turn.

\begin{proposition}[Non-degeneracy]
\label{p.change}
Let $\delta \in (0,1)$. When the bounds $M \geq c^{-2/3} \big( \log 4C \delta^{-1} \big)^{2/3}$  and $n \geq (M+1)^9 c^{-9} \vee c^{-2} \big( \log 4C \delta^{-1} \big)^2$ are satisfied,
$$
\PP \Big( Z_{n}(M)-Z_{n}(-M) \geq 4(2^{1/2} -  1) M \Big) \, \geq \, 1-\delta \, . 
$$ 
This assertion also holds when $Z_n$ is replaced by $Z$.
\end{proposition}

For $\e,K>0$ and $\alpha< 1/2$, a real-valued function $f$ whose domain contains $[-M,M]$ is said to be $(\e, K, \alpha)$-regular if, for all intervals $I\subset [-M,M]$ of length $\e$, $\sup_{y\in I} f(y)-\inf_{y\in I} f(y) \leq K\e^\alpha$. 

\begin{proposition}[Distribution of measure]\label{p.reg}
Let $\alpha \in (0, 1/2)$ and $M>0$. Almost surely, there exists a random value $\e_* > 0$ such that $Z$ is $(\e,8,\alpha)$-regular on $[-M,M]$ for all $\e\in (0,\e_*]$ .
\end{proposition}

{\bf Proof of Theorem~\ref{t.hd}(2).} 
We indeed take $E = [-M,M]$ and $\mu$ specified by $\mu(a,b] = Z(b) - Z(a)$ in Theorem~\ref{t.mdp}. Choosing $M \geq 4^{-1} (2^{1/2} - 1)^{-1}$ in Proposition~\ref{p.change}, and applying this result in the case of~$Z$, we see that $\mu(E) \geq 1$ with probability at least $1 - \delta$. (In fact, this lower bound of one is not needed; merely that $\mu(E) > 0$ would suffice.) From Proposition~\ref{p.reg}, we see that the hypothesis~(\ref{e.vcdelta}) is verified for any given $\alpha \in (0, 1/2)$ with $K = 8$ and for a random but positive choice of the constant $\eta$. Thus we find, as desired,  that the Hausdorff dimension of $\mathrm{LV}(Z) \cap [-M,M]$ is at least one-half with a probability that is at least $1 - \delta$. \qed

In order to prove Proposition~\ref{p.change}, we recall upper and lower tail bounds for the parabolically adjusted weight  $\wgt{n}{0}{z} + 2^{-1/2}  z^2$. The next result is quoted from~\cite{H16}, but it is a consequence of bounds on the upper and lower tails of the highest eigenvalue of a matrix randomly drawn from the Gaussian unitary ensemble, bounds respectively due to Aubrun~\cite{Aubrun} and Ledoux~\cite{Ledoux}. 

\begin{theorem}\cite[Proposition~$2.5$]{H16}
\label{t.onepoint}
 If $x,y \in \R$ satisfy $y-x \geq - 2^{-1} n^{1/3}$ and $\vert y - x \vert \leq c n^{1/9}$, then
$$
\PP \bigg( \Big\vert \wgt{n}{x}{y} + 2^{-1/2}  (y-x)^2 \Big\vert \geq s \bigg) \leq C \exp \big\{ - c s^{3/2} \big\}
$$
for all $s \in \big[1, n^{1/3} \big]$.
\end{theorem}
{\bf Proof of Proposition~\ref{p.change}.} The latter assertion of the proposition, concerning $Z$, follows from the former by the Skorokhod representation of weak convergence. To prove the former, consider
$x,y \in \R$ that satisfy $y - x \geq 2^{-1}n^{1/3}$,  and set 
$\omega(x,y) =  \wgt{n}{x}{y} + 2^{-1/2} (y-x)^2$ equal to the parabolically adjusted weight associated to the polymer $\polymer{n}{x}{y}$.  Note that 
$$
Z_n(M) = \omega(1 , M) - \omega(-1 , M) + 2^{3/2} M
 \, \, \, \, \textrm{and} \, \, \, \,  
 Z_n(-M) = \omega(1 , -M) - \omega(-1 , -M) - 2^{3/2} M \, .
 $$

Set $m_0 = c^{-2/3} \big( \log 4C \delta^{-1} \big)^{2/3}$.
Theorem~\ref{t.onepoint} implies that, when $m_0 \geq 1$, $M  > 0$ and  $n \geq (M+1)^9 c^{-9} \vee m_0^3$,
$$
\PP \Big( \max  \big\{ \vert \omega(1,M) \vert , \vert \omega(-1,M) \vert , \vert \omega(1,-M) \vert ,\vert \omega(-1,-M) \vert \big\} \geq   m_0 \Big) \leq \delta \, ,
$$
Suppose now that $M \geq m_0$. The four $\omega$ quantities are all at most $M$ in absolute value except on an event of probability at most $\delta$. In this circumstance, we have the bound $Z_n(M) - Z_n(-M) \geq (2^{5/2} -  4) M$, so that the proof of Proposition~\ref{p.change} is completed.  \qed

It remains only to validate our second tool, concerning distribution of measure.

{\bf Proof of Proposition~\ref{p.reg}.}
Let $Z_{\infty;\pm 1}$ denote a random function whose law is an arbitrary weak limit point of $z \to \wgt{n}{\pm 1}{z}$ as $n \to \infty$. 
First note that it suffices to prove that almost surely there exists $\e_*$ such that  $Z_{\infty;-1}$ and $Z_{\infty;1}$ are $(\e,4,\alpha)$-regular on $[-M,M]$ for all $\e\in (0, \e_*]$.  We will prove this for  $Z_{\infty;-1}$, the other argument being no different. With $\e_k = 2^{-k}$, it is moreover enough to argue that  there exists a random value $K_0 \in \N$ for which $Z_{\infty;-1}$ is $(\e_{k},2,\alpha)$-regular whenever $k\geq K_0$, since this implies that this random function is $(\e_{k},4,\alpha)$-regular for each $\e\leq \e_{K_0}$. 

Corollary \ref{c.aest}  implies that, for any $M$, there exists $k_0=k_0(M)$ such that, for $k\ge k_0$, 
\begin{equation}
\nonumber
 \limsup_{n \in \N} \, \P \left( \sup_{\begin{subarray}{c}y\in [-M,M] \, , \\
\eta_1,\eta_2 \in [0,2^{-k}]
\end{subarray}} 
       \Big\vert \wgt{n}{0}{y + \eta_2} - \wgt{n}{0}{y + \eta_1} \Big\vert  \geq 2^{-k\alpha}\right) \leq \exp \big\{ -2^{3k(1-2\alpha)/4} \big\} \, .
\end{equation}

Since this right-hand side is summable in $k$,
the Borel-Cantelli lemma implies that almost surely there exists a random positive integer $K_0$ such that on $[-M,M]$, the weak limit point $Z_{\infty;-1}$ is $(\e_{k},2,\alpha)$-regular for all $k\geq K_0$. This completes the proof of Proposition~\ref{p.reg}. \qed

\subsection{A lower bound on the probability of polymer pairs with close endpoints}
These last paragraphs are devoted to giving a remaining proof, that of Theorem~\ref{t.disjtpoly.lb}.
The derivation has three parts. First we state and prove Proposition~\ref{p.disjtpoly.lb}, which is an averaged version of the sought result. Then follows Proposition~\ref{p.polyscale}, which indicates all terms being averaged are about the same. From this we readily conclude that Theorem~\ref{t.disjtpoly.lb} holds.

To state our averaged  result, 
let $K$ and $\e$ be positive parameters; and write $\mc{I}(K,\e)$ for
the set of intervals of the form $[u,u+\e]$ that intersect $[-K,K]$ and for which $u \in \e \Z$.
The cardinality of $\mc{I}(1,\e) \times \mc{I}(K,\e)$ is of order $K \e^{-2}$, so that Proposition~\ref{p.disjtpoly.lb} indeed concerns the average value of 
the probability that $\maxpolymac{n}{I}{J}  \geq 2$ as $I$ and $J$ vary over those intervals in a compact region that abut consecutive elements of $\e \Z$.

\begin{proposition}\label{p.disjtpoly.lb} 
There exists $K_0 > 0$ such that, for $\eta > 0$ and $K \geq K_0$, we may find $\e_0 = \e_0(K,\eta)$
for which, whenever  $\e \in (0,\e_0)$, there exists 
$n_0 = n_0(K,\e,\eta)$ so that $n \in \N$, $n \geq n_0$, implies that
\begin{equation}\label{e.disjtpoly}
\e^2 \sum_{\begin{subarray}{c} 
     I \in \mathcal{I}(1,\e)\, , \\ J \in \mathcal{I}(K,\e)  \end{subarray}}\PP \Big( \maxpolymac{n}{I}{J}  \geq 2 \Big) 
  \, \geq  \,
2^{-3} \e^{3/2 + \eta} \, .
\end{equation}
\end{proposition}
{\bf Proof.}
The proposition asserts its result when $K \geq K_0$, $\e \in (0,\e_0)$ and $n \in \N$ satisfies $n \geq n_0$. We begin the proof by noting that explicit choices of these three bounds on parameters will be seen to be given by 
$K_0 =  c^{-2/3} \big( \log 8C  \big)^{2/3}$;  
$$
\e_0  = 2^{-1} \min \Big\{ 10^{-13} K^{-2}  C^{-2}  \big( c_1 2^{-19} 2^{3/2} \eta^{-1} \big)^{2/\eta} , \big( 2^{1/2}(K+1) \big)^{-1/(1 + \eta)} ,  \big( 5 \cdot 10^3 \big)^{-1/\eta}\Big\} \, ;
$$
and
 $n_0  = \max \big\{ 2^{36} 3^{18} c^{-18}(K+1)^{18}, 2^{18} 10^{-6} \e^{-18\eta}, (K+1)^9 c^{-9} , c^{-2} \big( \log 8C  \big)^2 \big\}$. 
 
Note first that 
$$
Z_n(K) - Z_n(-K) \, = \, \Big( W_n(1,K) - W_n(1,-K) \Big) - \Big( W_n(-1,K) - W_n(-1,-K) \Big)
$$
is at most 
$$
  \sum_{J \in \mathcal{I}(K,\e)} \, \sup_{\begin{subarray}{c} 
    v_1,v_2 \in J  \end{subarray}} \, \Big( \, \big\vert \wgt{n}{-1}{v_1}  - \wgt{n}{-1}{v_2}  \big\vert  + \big\vert \wgt{n}{1}{v_1}  - \wgt{n}{1}{v_2}  \big\vert \, \Big) \cdot {\bf 1}_{J \cap \mathrm{LV}(Z_n)  \not= \emptyset} \, ,  
$$
where the indicator function ${\bf 1}$ may be included because $Z_n(K) - Z_n(-K)$
may be viewed as a telescoping sum of differences indexed by intervals $J \in \mathcal{I}(K,\e)$ 
of which those disjoint from $\mathrm{LV}(Z_n)$ contribute zero. 

Let $J \in  \mathcal{I}(K,\e)$ be given.
We now apply Theorem~\ref{t.differenceweight} 
with parameter choices $x = -1$; $y \in [-K ,K ]$ the left endpoint of $J$; and $R = 2\e^{-\eta}$.
By supposing that $\e \leq \big( 2^{1/2}(K+1) \big)^{-1/(1 + \eta)}$,
the parabolic term $\big\vert Q(v-u) - Q(y-x) \big\vert$ in this theorem is at most $2^{1/2} (K+1) \e \leq \e^{-\eta}$, so that the theorem implies that
$$
\PP \Big( \sup_{\begin{subarray}{c} 
    v_1,v_2 \in J  \end{subarray}} \,  \big\vert \wgt{n}{-1}{v_1}  - \wgt{n}{-1}{v_2}  \big\vert  \geq \e^{1/2 - \eta} \Big) \, \leq \,  10032 \, C \exp \big\{ - c_1 2^{-19} 2^{3/2} \e^{-3\eta/2} \big\}
$$
 provided that $n \geq  \max \big\{ 10^{32}c^{-18} , 2^{36} 3^{18} c^{-18} (K+1)^{18}, 2^{18} 10^{-6} \e^{-18\eta} \big\}$ and $\e \leq \big( 5 \cdot 10^3 \big)^{-1/\eta}$. We may equally apply Theorem~\ref{t.differenceweight}  with $x=1$ to find that the same estimate holds when the quantity $\big\vert \wgt{n}{1}{v_1}  - \wgt{n}{1}{v_2}  \big\vert$ is instead considered.
 
 Setting $\mathsf{G}$ to be the event that  $\sup_{\begin{subarray}{c} 
    v_1,v_2 \in J  \end{subarray}} \,  \big\vert \wgt{n}{x}{v_1}  - \wgt{n}{x}{v_2}  \big\vert  < \e^{1/2 - \eta}$ holds for all $x \in \{-1,1\}$ and $J \in \mathcal{I}(K,\e)$, we see that, on $\mathsf{G}$, 
    $$
Z_n(K) - Z_n(-K) \, \leq \, 2 \e^{1/2 - \eta} \cdot
  \Big\vert \Big\{ J \in \mathcal{I}(K,\e) J \cap \mathrm{LV}(Z_n)  \not= \emptyset \Big\} \Big\vert \, , 
    $$
    and that 
    \begin{equation}\label{e.goodcomp}
     \PP \big( \mathsf{G}^c \big) \leq \big( 2K \e^{-1} + 1  \big) \cdot  2 \cdot 10032 \, C \exp \big\{ - c_1 2^{-19} 2^{3/2} \e^{-3\eta/2} \big\} \, .
    \end{equation}

By taking $M = K \geq 4^{-1} (2^{1/2} - 1)^{-1}$ in Proposition~\ref{p.change}, our choice of $K \geq c^{-2/3} \big( \log 8C  \big)^{2/3}$ ensures that,
when $n\in \N$ satisfies $n \geq (M+1)^9 c^{-9} \vee c^{-2} \big( \log 8C  \big)^2$, it is with probability at least one-half that  the event $Z_n(K)-Z_n(-K) \geq 1$ occurs.    
We thus find that 
$$
 \PP \bigg( 2 \e^{1/2 - \eta} \cdot
  \Big\vert \Big\{ J \in \mathcal{I}(K,\e): J \cap \mathrm{LV}(Z_n)  \not= \emptyset \Big\} \Big\vert \geq 1 \bigg) \geq 1/4 \, ,
$$
provided that the right-hand side of~(\ref{e.goodcomp}) is at most $1/4$ -- as it is, due to a brief omitted estimate that uses 
$\e \leq 1$, $K \geq 1$ and the hypothesised upper bound $\e \leq 10^{-13} K^{-2}  C^{-2}  \big( c_1 2^{-19} 2^{3/2} \eta^{-1} \big)^{2/\eta}$.
We see then that 
$$
\E \, \Big\vert \Big\{ J \in \mathcal{I}(K,\e): J \cap \mathrm{LV}(Z_n)  \not= \emptyset \Big\} \Big\vert \, \geq \, 2^{-3} \e^{-1/2 + \eta} \, . 
$$

Proposition~\ref{p.mdp} implies that, for $J\in \mathcal{I}(K,\e)$,
$$
\PP \Big( J \cap \mathrm{LV}(Z_n)  \not= \emptyset \Big) \, \leq \, \sum_{I \in \mathcal{I}(1,2\e)}\PP \Big( \maxpolymac{n}{I}{J}  \geq 2 \Big) \, .
$$
Thus,
$$
 \sum_{\begin{subarray}{c} 
     I \in \mathcal{I}(1,2\e)\, , \\ J \in \mathcal{I}(K,\e)  \end{subarray}}\PP \Big( \maxpolymac{n}{I}{J}  \geq 2 \Big) \geq  2^{-3} \e^{-1/2 + \eta} \, .
$$
The conclusion of Proposition~\ref{p.disjtpoly.lb} would be achieved were $\mathcal{I}(1,2\e)$ to read $\mathcal{I}(1,\e)$. We relabel $\e$ to be the present $2\e$ in order to achieve this; note that it is this relabelling which is responsible for the presence of a factor of one-half in the specification of the value of $\e_0$ at the beginning of the proof. \qed

The next result -- that the terms being averaged are all roughly equal -- is inspired by the first line of page $34$ of the second version of~\cite{DOV18}.
\begin{proposition}\label{p.polyscale}
Let $n \in \N$ and $x,y \in \R$ satisfy $y - x \geq - 2^{-1}n^{1/3}$. Then
$$
\PP \Big( \maxpolymac{n}{[x,x+\e]}{[y,y+\e]}  \geq 2 \Big) \, = \,
\PP \Big( \maxpolymac{n}{[0,\e']}{[0,\e']}  \geq 2 \Big) \, ,
$$
where $\e' > 0$ is a quantity that differs from $\e$
by at most $O(1) n^{-1/3} \big( \vert y-x \vert + \e \big) \vert y -x \vert$. The constant factor implied by the use of the $O(1)$ notation has no dependence on $(n,x,y)$.
\end{proposition}
\noindent{\bf Proof.} Let $n \in \N$ and $x \in \R$. Since the Brownian paths in the underlying environment $B:\Z \times \R \to \R$ have stationary increments, we may replace this ensemble by the system $\Z \times \R \to \R: (k,z) \to B(k,z - 2n^{2/3}x)$ without changing the ensemble's law. By the form of the scaling map $R_n:\R^2 \to \R^2$, we find that
$$
\PP \Big( \maxpolymac{n}{[x,x+\e]}{[y,y+\e]}  \geq 2 \Big) = 
\PP \Big( \maxpolymac{n}{[0,\e]}{[y-x,y-x+\e]}  \geq 2 \Big) \, .
$$
It is thus enough to prove Proposition~\ref{p.polyscale} in the case that $x = 0$. To this end, we let $n \in \N$ and $y \in  \R$ be given. Consider again the ensemble $B:\Z \times \R \to \R$. Set $B'(k,z) = B(k,\alpha z)$ with $\alpha = 1 + 2 n^{-1/3}y$. 
Equally, we may write $B' = \alpha^{1/2}B$, and this identity shows us that LPP under $B$ and $B'$ differ merely by a multiplication of energy by a factor of $\alpha^{1/2}$; so that the change $B \to B'$ makes no difference in law to the geometry of geodesics including their disjointness.

A geodesic specified by the ensemble $B':\Z \times \R \to \R$ that runs between $(0,0)$ and $(n,n)$ corresponds to a geodesic specified by  $B:\Z \times \R \to \R$
that runs between $(0,0)$ and $(n + 2n^{2/3}y,n)$. When the scaling map $R_n$ is applied, polymers $(0,0) \to (y,1)$ and $(0,0) \to (0,1)$ result from the primed and original environments. On the other hand, an original geodesic running between  $(2n^{2/3}\e,0)$ and $\big(n + 2n^{2/3}(y+\e),n\big)$ corresponds to a primed geodesic between $(2n^{2/3}\e,0)$ and $\big(n,n+2n^{2/3}\e + \gamma\big)$, where the small error $\gamma$ is readily verified to satisfy $\gamma = O(1) n^{1/3}(\vert y \vert + \e)\vert y \vert$. Applying the scaling map again, original and primed polymers are seen to run respectively $(\e,0) \to (\e,1)$ and $(\e',0) \to (\e',1)$, where $\e' > 0$ satisfies $\vert \e - \e' \vert = O(1) n^{-1/3} (\vert y \vert + \e)\vert y \vert$. Thus do we confirm Proposition~\ref{p.polyscale} in the desired special case that $x = 0$. \qed


\noindent{\bf Proof of Theorem~\ref{t.disjtpoly.lb}.}
By Proposition~\ref{p.polyscale}, $\PP \Big( \maxpolymac{n}{[0,2\e]}{[0,2\e]}  \geq 2 \Big)$
is at least the value of every summand on the left-hand side of~(\ref{e.disjtpoly}), provided that $n$ exceeds a level determined by $K$ and $\e$. Since $\vert \mc{I}(K,\e) \vert \leq 2K\e^{-1} + 2$, we see that Proposition~\ref{p.disjtpoly.lb} with $K = K_0$ implies Theorem~\ref{t.disjtpoly.lb} with 
$d = 
2^{-7} K_0^{-1}$.   \qed

\bibliography{airydifference}
\bibliographystyle{plain}

\end{document}